\documentclass[10pt]{article}

\usepackage{amsmath}
\usepackage{mathtools}
\usepackage{amssymb}
\usepackage{latexsym}
\usepackage{enumitem}
\usepackage{mathrsfs}
\usepackage{comment}
\usepackage{pifont}
\usepackage{enumitem}
\usepackage{bbold}
\usepackage{graphicx}
\usepackage{color}
\usepackage{hyperref}
\hypersetup{
    colorlinks=true,
    linkcolor=blue,
    filecolor=magenta,      
    urlcolor=cyan,
    citecolor=magenta
}

\newtheorem{assumption}{Assumption}

\def\qed{ \ \vrule width.2cm height.2cm depth0cm\smallskip}
\newenvironment{proof}{\noindent {\bf Proof.\/}}{$\qed$\vskip 0.1in}

\newcommand{\la}{\langle}
\newcommand{\ra}{\rangle}

\newcommand{\ol}{\overline}

\newcommand{\ba}{\begin{array}}
\newcommand{\ea}{\end{array}}
\newcommand{\be}{\begin{equation}}
\newcommand{\ee}{\end{equation}}
\newcommand{\bea}{\begin{eqnarray}}
\newcommand{\eea}{\end{eqnarray}}
\newcommand{\beaa}{\begin{eqnarray*}}
\newcommand{\eeaa}{\end{eqnarray*}}

\def\dbE{\mathbb{E}}
\def\dbF{\mathbb{F}}

\def\dbI{\mathbb{I}}

\def\dbL{\mathbb{L}}

\def\dbP{\mathbb{P}}
\def\dbR{\mathbb{R}}

\def\dbT{\mathbb{T}}

\def\a{\alpha}
\def\b{\beta}
\def\g{\gamma}
\def\d{\delta}
\def\e{\varepsilon}

\def\k{\kappa}
\def\l{\lambda}
\def\m{\mu}
\def\n{\nu}
\def\si{\sigma}
\def\t{\tau}
\def\f{\varphi}
\def\th{\theta}

\def\G{\Gamma}
\def\D{\Delta}

\def\O{\Omega}
\def\cB{{\cal B}}

\def\cF{{\cal F}}

\def\cJ{{\cal J}}

\def\cL{{\cal L}}

\def\cN{{\cal N}}

\def\cP{{\cal P}}

\def\cS{{\cal S}}

\def\cW{{\cal W}}

\def\no{\noindent}

\def\q{\quad}
\def\qq{\qquad}

\def\pa{\partial}
\def\cd{\cdot}
\def\cds{\cdots}

\def\tr{\hbox{\rm tr}}

\def\qed{ \hfill \vrule width.25cm height.25cm depth0cm\smallskip}

\newcommand{\basa}{\begin{assumption}}
\newcommand{\easa}{\end{assumption}}

\newcommand{\bas}{\begin{assum}}
\newcommand{\eas}{\end{assum}}

\def\pa{\partial}

 \def\cd{\cdot}
\def\cds{\cdots}

\def\supp{\hbox{\rm supp$\,$}}

\def\tr{\hbox{\rm tr$\,$}}

\def\dis{\displaystyle}

\def\bQ{{\bf Q}}
\def\bx{{\bf x}}
\def\cad{c\`{a}dl\`{a}g}

\def\1{{\bf 1}}

\def\:{\!:\!}
\def\reff#1{{\rm(\ref{#1})}}
\def \proof{{\noindent \bf Proof\quad}}

\def \bS{{\bf S}}
\def \bm{{\bf m}}

at 9pt

\begin{document}

\newtheorem{thm}{Theorem}[section]
\newtheorem{lem}[thm]{Lemma}
\newtheorem{cor}[thm]{Corollary}
\newtheorem{prop}[thm]{Proposition}
\newtheorem{rem}[thm]{Remark}
\newtheorem{eg}[thm]{Example}
\newtheorem{defn}[thm]{Definition}
\newtheorem{assum}[thm]{Assumption}

\numberwithin{equation}{section}

\title{\bf{Dynamic programming equation for the mean field optimal stopping problem\footnote{The first two authors are grateful for the financial support from the Chaires FiME-FDD and Financial Risks of the Louis Bachelier Institute. The third author is supported in part by NSF grant DMS-1908665.}}}
\author{Mehdi Talbi\footnote{Department of Mathematics, ETH Zürich, Switzerland, mehdi.talbi@math.ethz.ch} \quad Nizar Touzi\footnote{CMAP, \'{E}cole polytechnique, France, nizar.touzi@polytechnique.edu}  \quad Jianfeng Zhang\footnote{Department of Mathematics, University of Southern California, United States, jianfenz@usc.edu. }}
\date{\today}

\maketitle

\begin{abstract}
We study the optimal stopping problem of McKean-Vlasov diffusions when the criterion is a function of the law of the stopped process. A remarkable new feature in this setting is that the stopping time also impacts the dynamics of the stopped process through the dependence of the coefficients on the law. The mean field stopping problem is introduced in weak formulation in terms of the joint marginal law of the stopped underlying process and the survival process. This specification satisfies a dynamic programming principle. The corresponding dynamic programming equation is an obstacle problem on the Wasserstein space, and is obtained by means of a general It\^o formula for flows of marginal laws of c\`adl\`ag semimartingales. Our verification result characterizes the nature of optimal stopping policies, highlighting the crucial need to randomized stopping. The effectiveness of our dynamic programming equation is illustrated by various examples including the mean-variance optimal stopping problem.
\end{abstract}

\no{\bf MSC2020.} 60G40, 49N80, 35Q89, 60H30

\vspace{3mm}
\no{\bf Keywords.} Mean field optimal stopping, McKean-Vlasov SDEs, dynamic programming.

\section{Introduction}

In this paper we study a McKean-Vlasov type of optimal stopping problem, where the state dynamics and/or the reward function depend on the law of the stopped process. To be precise, given $X_0$ and an independent Brownian motion $W$, consider
\bea\label{mkvintro}
X_t = X_0 + \int_0^{t \wedge \t} b(s,X_s,\cL_{X_s})ds + \int_0^{t \wedge \t} \si(s,X_s,\cL_{X_s})dW_s,
\eea
where $\t$ is a stopping time and $\cL_{X_s}$ denotes the law of $X_s$. We emphasize the impact of $\t$ on $\cL_{X_s}$, in particular,  $\cL_{X_s}$ is neither equal to $\cL_{X^0_{\t\wedge s}}$  nor to $\cL_{X^0_s}|_{s=\t}$, where $X^0$ denotes the unstopped process:
\bea
\label{X0}
X^0_t = X_0 + \int_0^t  b(s,X^0_s,\cL_{X^0_s})ds + \int_0^t  \si(s,X^0_s,\cL_{X^0_s})dW_s.
\eea
Our optimization problem is, for some functionals $f$ and $g$ defined on a space of probability laws,
 \bea\label{generic}
V_0 := \sup_{\t} \dbE\Big[ \int_0^\t f(s, X_s, \cL_{X_s})ds \Big] + g(\cL_{X_\t}).
\eea
When $b, \si$ and $f$ do not depend on $\cL_{X_s}$ and $g(\cL_{X_\t}) = \dbE\big[\f(X_\t)\big]$ for some function $\f: \dbR^d\to \dbR$, the above problem reduces to a standard optimal stopping problem, see e.g. Shiryaev \cite{Shiryaev}.  The mean field optimal stopping problem \reff{generic} can be viewed as the limit of a multiple stopping problem over a large system interacting through the empirical measure:  
\bea\label{Nparticles}
\left.\ba{c}
\dis X_t^i = x_i + \int_0^{t \wedge \t_i} b(s, X_s^i, \bar \m_s)ds + \si(s, X_s^i, \bar \m_s)dW_s^i, \q \bar \mu_s:= {1\over N}\sum_{i=1}^N \d_{X^i_s};\\
\dis V_0^N := \sup_{(\t_1, \dots, \t_N)} \dbE\Big[\frac 1 N \sum_{i=1}^N \int_0^{\t_i} f(s, X_s^i,\bar \m_s)ds + g\Big(\frac 1 N \sum_{i=1}^N \d_{X_{\t_i}^i} \Big)\Big],
\ea\right.
\eea
where $\d_x$ denotes the Dirac-measure, $(W^1, \dots, W^N)$ are $N \times d-$dimensional Brownian motions. We refer to Kobylanski, Quenez \& Rouy-Mironescu \cite{KQR} for general multiple stopping problems, and we shall investigate the convergence issue in an accompanying paper \cite{TTZ3}. 

There has been a strong attention on mean field games of optimal stopping in the literature, see, e.g., Bertucci \cite{ber}, Bouveret, Dumitrescu \& Tankov \cite{BDT}, Carmona, Delarue \& Lacker \cite{CDL}, and Nutz \cite{nutz}. Given $\{\mu_t\}_{t\ge 0}$, consider the optimal stopping problem:
\bea
\label{Vmu0}
V^{\mu_\cd}_0:= \sup_{\t}  \dbE\Big[ \int_0^\t f(s, X_s^{\m_.}, \m_s)ds + g(\t, X^{\mu_\cd}_\t, \mu_\t) \Big], % + g(\cL_{X^{\mu_\cd}_\t}), 
\eea
where  $X^\mu$ is  unstopped and solves a standard SDE (not McKean-Vlasov type as in \reff{X0}): 
\beaa
X^{\mu_\cd}_t = X_0 + \int_0^t b(s,X^{\mu_\cd}_s,\mu_s)ds + \int_0^t \si(s,X^{\mu_\cd}_s,\mu_s)dW_s.
\eeaa
Assume the above problem has an optimal stopping time $\t^*(\mu_\cd)$, then the mean field game problem is to find a fixed point $\{\mu_t\}_{t\ge 0}$, namely the mean field equilibrium: $\cL_{X^{\mu_\cd}_{\t^*(\mu_\cd)\wedge t}} = \mu_t$, $t\ge 0$. We remark that in the last mean field game,  for given $\{\mu_t\}_{t\ge 0}$, the dynamics of $X^{\mu_\cd}$ does not depend on the stopping time $\t$  and the optimal stopping problem \reff{Vmu0} is a standard one as in \cite{Shiryaev},  so it has a completely different structure than our optimal stopping problem. We would also like to mention  Li \cite{Li1},  Briand, Elie \& Hu \cite{BEH}, and Djehiche, Elie \& Hamadene \cite{DEH} for closely related works on mean field type reflected BSDEs, and Belomestny \& Schoenmakers \cite{BS} for a numerical method for mean field type optimal stopping problems. However, in all these works again the dynamics of the state process  does not depend on the stopping time $\t$.  To our best knowledge, our work is the first in the literature to study the optimal stopping problem where the dynamics depends on the law of the stopped process, or say in \reff{Nparticles} the interaction is through the stopped particles. 

Besides the obvious connection with large interacting particle systems, the general form \reff{generic} is convenient for many other applications. For example, by considering the unstopped state process $X^0$ in \reff{X0},  the optimal stopping of mean variance problem $\sup_{\t} \big\{\dbE[X_\t^{0}]-\frac{1}{2}\mathrm{Var}(X_\t^{0})\big\}$ corresponds to $g(\mu) = \int_{\dbR} (x -{1\over 2} x^2)\mu(dx) + \big(\int_{\dbR} x \mu(dx)\big)^2$ for a square integrable measure $\m$. Another example is the optimal stopping problem under probability distortion, used in behavioral economics, which corresponds to $g(\mu) = \int_0^\infty \f\big(\mu([U^{-1}(y), \infty)\big)dy,$ for some utility function $U:\dbR\to [0, \infty)$, and some distortion function $\f : [0,1] \longrightarrow [0,1]$. When $X^0$ is a Geometric Brownian motion and the time horizon is infinite, Pedersen \& Peskir \cite{PP} proved the existence of optimal stopping time for the mean variance problem,  and Xu \& Zhou \cite{XZ} obtained the optimal stopping time for the probability distortion problem for some special shapes of the functions $\f$ and $U$ (convex, concave, or reverse S-shaped). We remark that these problems are typically considered as time inconsistent problems, as we will explain in the next paragraph, and the existing literature considers only the static problem, namely the existence of optimal stopping time for the problem over a fixed time interval ($[0, \infty)$ or $[0, T]$). We shall study the problem \reff{generic} systematically, and more importantly, dynamically. We remark that, even when we consider only the unstopped state process $X^0$, our dynamic approach for the optimal stopping problem \reff{generic} seems new.
  
It is well known that standard optimal stopping problems can be solved by the dynamic programming approach, see e.g. El Karoui \cite{EK}  and Shiryaev 
\cite{Shiryaev}. The situation here is more subtle because of the involvement of the law. In order to have Dynamic Programming Principle (DPP, for short), it is crucial to choose the right variable, which stands for the information one needs  to make the dynamic system ‘‘Markovian". Indeed, if we define $V(t,x)$ as the dynamic value function for problem \reff{generic} on $[t, T]$ with initial condition $X_t=x$, which in the case \reff{Nparticles} means we observe only the state $x_i$ of one particular player $i$, the DPP would fail.  Consequently the problem is often viewed as time inconsistent in the standard sense. Moreover, even if we define $V(t,\mu)$ as the dynamic value function for problem \reff{generic} on $[t, T]$ with initial condition $\cL_{X_t}=\mu$, the DPP would still fail. 

Our first observation is that a successful DPP requires the introduction of the survival process $I_t := \1_{\{\t>t\}}$. To be precise, we will have the desired DPP if we write the dynamic value function as $V(t, \cL_{(X_t, I_t)})$, that is, to maintain the time consistency, we need to know not only the current states of all particles, but also which particles are still surviving. Moreover, we formulate a weak relaxed version of \eqref{mkvintro} by allowing for randomized stopping times induced by the set $\cP(t,m)$ of all joint distributions $\dbP$ of the stopped process and the corresponding stopping time, started at time $t$ from the initial distribution $m$. Such a weak formulation is particularly convenient here for two reasons: 

\no$\bullet$ the set of controls has been shifted from the stopping times into $\cP(t,m)$, that we will prove to  be compact, implying the existence of an optimal $\dbP^*$ to the mean field optimal stopping problem as long as $f$ and $g$ are upper-semicontinuous; 

\no$\bullet$ shifting the state variable from the process $X$ into the flow of joint marginal distributions, denoted as  $\{\dbP_{(X_t,I_t)}\}$ in order to emphasize its dependence on $\dbP$, enables us to establish a DPP and to derive a dynamic programming equation on the space of measures to characterize  the value function $V$.

\no More precisely, given that the laws are deterministic, our following DPP is very easy to establish:
$$
V(t, m) = \sup_{\dbP\in \cP(t, m)}\int_t^s \dbE[f(r, X_r, \dbP_{(X_r, I_r)}) I_r] dr+ V\big(s, \dbP_{(X_s, I_s)}\big),
$$ 
Such dynamic programming approach has also been used successfully in the mean field control literature, where the state variable is $\cL_{X_t}$, see, e.g.,   Carmona \& Delarue \cite[Vol. 1, ch. 6]{CarDel}, Pham \& Wei \cite{PW}, Wu \& Zhang \cite{WZ}, and Djete, Possamai \& Tan \cite{DPT}. 

The corresponding dynamic programming equation is as usual derived by means of Itô's formula. Itô's formula for functions on Wasserstein space of probability measures has been established for continuous diffusions by Buckdahn, Li, Peng \& Rainer \cite{BLPR} and Chassagneux, Crisan \& Delarue \cite{CCD}, and for jump diffusions by Li \cite{Li2} and Burzoni, Ignazio, Reppen \& Soner \cite{BIRS}. However, \cite{Li2, BIRS} require the law of the state process to be continuous under the Wasserstein distance, while in our case it is quite possible that $t\longmapsto \dbP_{I_t}$ is discontinuous. We thus first extend Itô's formula so that both the state process and its law can have jumps. Our proof follows the standard derivation, based on the linear functional derivative. We introduce an appropriate time discretization and reduce our derivation to the standard Itô's formula for càdlàg semimartingales. We also refer to the independent work of Guo, Pham \& Wei \cite{GPW}, who prove similar results by using density arguments under slightly different technical conditions, see Remark \ref{rem:gpw}.

Together with the DPP, our Itô's formula immediately leads to the desired dynamic programming equation, an obstacle problem on the Wasserstein space. We shall characterize the value function, provided its sufficient regularity, as the unique classical solution of the obstacle problem, and we will use the value function to characterize the structure of the optimal stopping time. The regularity of the value function, of course, remains a challenging problem in general, and we will therefore investigate the viscosity solution approach for the obstacle problem in another accompanying paper \cite{TTZ2}. 
 
The paper is structured as follows. In Section \ref{meanfield}, we set the mean field optimal stopping problem in weak formulation, and establish the dynamic programming principle. In Section \ref{Ito} we prove the It\^o's formula for possibly discontinuous flows of measures of semimartingales, that in particular allows us to differentiate smooth functions along the flow $\{\dbP_{(X_t,I_t)}\}_{t \in [0,T]}$. In Section \ref{obs} we derive the dynamic programming equation for the value function and establish its classical solution theory. 
Section \ref{examples} is dedicated to some examples illustrating the connection with the standard optimal stopping theory, and shedding more light on a class of criteria including the mean-variance one. We also provide an explicit example which exhibits both features of pure stopping strategies and randomized ones. In Section \ref{sect-extension} we provide two extensions.  Subsection \ref{infinite} extends our results to the infinite horizon setting, and Subsection \ref{JumpDiffSec} provides a quick discussion of the extension  to the case where the process $X$ is a jump-diffusion. Finally, Appendices \ref{appendixA} and \ref{appendixB} report some technical proofs. 

\vspace{2mm}
\no {\bf Notations.} We denote by $\cP(\O,\cF)$ the set of probability measures on a measurable  space $(\O,\cF)$, and $\cP_2(\O,\cF)$ the subset of square integrable probability measures in $\cP(\O,\cF)$, equipped with the $2$-Wasserstein distance $\cW_2$. When $(\O,\cF) = (\dbR^d,\cB(\dbR^d))$, we simply denote them as $\cP(\dbR^d)$ and $\cP_2(\dbR^d)$. For a random variable $Z$ and a probability $\dbP$, we denote by $\dbP_Z:=\dbP\circ Z^{-1}$ the law of $Z$ under $\dbP$. 
 For vectors $x, y\in \dbR^n$ and matrices $A, B\in \dbR^{n\times m}$, denote $ x\cd y:=\sum_{i=1}^n x_iy_i$  and $A:B:= \tr(A B^\top)$.

\section{Formulation of the mean field optimal stopping problem}\label{meanfield}

Let $T < \infty$ be fixed, and $\O := C^0([-1,T],\dbR^d) \times \dbI^0([-1,T])$ the canonical space, where: 

\no$\bullet$ $C^0([-1,T],\dbR^d)$ is the set of continuous paths from $[-1,T]$ to $\dbR^d$, constant on $[-1,0)$; \\
\no$\bullet$ $\dbI^0([-1,T])$ is the set of non-increasing and càdlàg maps from $[-1,T]$ to $\{0,1\}$,  constant  on $[-1,0)$, and ending with value $0$ at $T$.
 
 \no We equip $\O$ with the Skorokhod distance, under which it is a Polish space. The choice of the extension to $-1$ is arbitrary,  the extension of time  to the left of the origin is only needed to allow for an immediate stop at time $t=0$.

We denote $Y:=(X,I)$ the canonical process, with state space $\mathbf{S}:=\dbR^d\times\{0,1\}$, its canonical filtration $\dbF = (\cF_t)_{t \in [-1,T]}$, and the corresponding jump time of the survival process $I$:
\bea
\label{tau}
\t := \inf\{t \ge 0 : I_t = 0\},  \ \mbox{so that $I_t := I_{0-}\1_{t < \t}$ for all $t \in [-1,T]$.}
\eea 
By the càdlàg property of $I$, $\t$ is an $\dbF-$stopping time. Denote further 
 \beaa
 \mathbf{Q}_t:=[t,T)\times\cP_2(\mathbf{S}), 
 &\mbox{and}&
 \overline{\mathbf{Q}}_t:=[t,T]\times\cP_2(\mathbf{S}),
 ~~t\in[0,T).
 \eeaa

 Let $(b,\si,f): [0,T] \times \dbR^d \times \cP_2(\mathbf{S}) \rightarrow \dbR^d \times \cS_d^+ \times \dbR$ and $g : \cP_2(\dbR^d) \rightarrow \dbR$, where  $\cS_d^+$ denotes the set of $d\times d$  non-negative symmetric matrices.  Throughout the paper, the following assumption will always be in force,  where $\cP_2(\bS)$ is equipped with the $\cW_2$-distance.  
\begin{assum}
\label{assum-bsig}
\no{\rm (i)} $b, \si$ are continuous in $t$, and uniformly Lipschitz continuous in $(x, m)$. 

{\rm (ii)}  $f$ is Borel measurable and has quadratic growth in $x\in \dbR^d$, and the following function $F$ is continuous on $[0, T]\times \cP_2(\bS)$: 
\bea
\label{F}
F(t,m) := \int_{\dbR^d} f(t,x,m)m(dx,1).
\eea
{\rm (iii)}  $g$ is upper-semicontinuous and locally bounded; and extended to $\cP_2(\bS)$ by $g(m) := g(m(\cd, \{0,1\}))$.
\end{assum}

Define the stopped McKean-Vlasov dynamics on $[0, T]$:
\bea
\label{WP}
X_s = X_0 + \int_0^s b(r, X_r, \dbP_{Y_r})I_rdr + \int_0^s \si(r, X_r, \dbP_{Y_r})I_rdW_r^\dbP \ \mbox{and} \ I_s = I_{0-}\1_{s < \t},
\eea
where a solution $\dbP$ of the last SDE is defined by the requirement that the following processes $M$ and $N$ are $\dbP-$martingales on $[0, T]$:
 \bea\label{martingalepb}
M_. := X_. - \int_0^. b(r,X_r,\dbP_{Y_r})I_rdr \ \mbox{and} \ N_.:= M_.^2 -  \int_0^. \si^2(r,X_r,\dbP_{Y_r})I_rdr.
\eea
 Note that $X_. = X_{. \wedge \t}$, and in particular $X_T = X_\t,$ $\dbP-$a.s. 
 
 We then focus on the mean field optimal stopping problem: given $\m \in \cP_2(\dbR^d)$,
\bea\label{naive}
V_0 := \sup_{\dbP} \dbE^\dbP\Big[\int_0^\t  f(r, X_r, \dbP_{Y_r})dr \Big] + g(\dbP_{X_\t}) = \sup_{\dbP} \int_0^T F(r, \dbP_{Y_r})dr + g(\dbP_{Y_T}) ,
 \eea
where the supremum is taken over all solutions $\dbP$ of the McKean-Vlasov SDE satisfying the constraint $\dbP_{X_0} = \mu$ and $\dbP(I_{0-}=1)=1$. We recall that this problem is motivated by the $N$-multiple optimal stopping problem \eqref{Nparticles}, whose convergence is studied in our accompanying paper \cite{TTZ3}. 
 
  In order to solve this problem, we use the dynamic programming approach, made possible by an appropriate dynamic version of the problem. This requires to take as a state the joint distribution $m_t$ of the variables $Y_t=(X_t,I_t)$, which leads
to the dynamic value function
 \begin{equation}\label{weakoptstop}
 V(t,m) := \underset{\dbP \in \cP(t,m)}{\sup} \int_t^T F(r, \dbP_{Y_r})dr  + g(\dbP_{Y_T}), \quad \mbox{$(t,m) \in \ol \bQ_0$},
 \end{equation}
where $\cP(t,m)$ is the set of probability measures $\dbP$ on $(\O,\cF_T)$ such that 

\no$\bullet$ $\dbP_{Y_{t-}} = m$ and $s\in [-1, t) \to Y_s$ is constant, $\dbP$-a.s.

\no$\bullet$ The processes $M,N$ of \eqref{martingalepb} are $\dbP$-martingales on $[t,T]$, so that, for some $\dbP$-Brownian motion $W^\dbP$,
\begin{equation}\label{asympt}
X_s = X_t + \int_t^s b(r, X_r, \dbP_{Y_r})I_r dr + \sigma(r, X_r, \dbP_{Y_r})I_r dW_r^\dbP , \ I_s = I_{t-} \1_{s < \t}, \ \dbP-\mbox{a.s.}
\end{equation}

\begin{prop}\label{existence} For any $(t,m) \in  \bQ_0$, the set $\cP(t,m)$ is compact under the Wasserstein distance $\cW_2$. Consequently, existence holds for the mean field optimal stopping problem \eqref{weakoptstop}.
\end{prop}
 We relegate this proof to Appendix \ref{appendixA}. Our main result of this section is the following dynamic programming principle (DPP for short).
\begin{thm}
\label{thm-DPP}
For any $(t,m)\in \bQ_0$ and $s\in [t, T]$, we have the DPP:
\bea\label{weakDPP}
 V(t,m) = \sup_{\dbP \in \cP(t,m)}  \int_t^s F(r, \dbP_{Y_r})dr  + V(s, {\dbP}_{Y_{s-}}) = \sup_{\dbP \in \cP(t,m)}  \int_t^s F(r, \dbP_{Y_r})dr  + V(s, {\dbP}_{Y_{s}}).
\eea 
\end{thm}

  \proof    Denote, for any probability measure $\dbP$ on $(\O, \cF_T)$,
  \beaa
  J(t, \dbP) := \int_t^T F(r, \dbP_{Y_r})dr  + g(\dbP_{Y_T}).
  \eeaa
 We start with proving the first equality of  \reff{weakDPP}.  Let $\tilde V(t, m)$ denote the middle term of \reff{weakDPP}. Fix an arbitrary $\dbP\in \cP(t, m)$, and denote $\tilde m:= \dbP_{Y_{s-}}$.
  
  First, for any time partition $\pi:-1=t_0<\cds<t_m < s \le t_{m+1}<\cds<t_{m+n} = T$, introduce the finite measure: for any $A_i\in \cB(\bS)$,
\beaa
\nu_\pi(A_0\times \cds \times A_{m+n}) := \dbP\Big(Y_{s-} \in \cap_{i=0}^{m} A_i, ~ Y_{t_{m+j}} \in A_{m+j}, j=1,\cds, n\Big). 
\eeaa
It is clear that $\{\nu_\pi\}_\pi$ satisfies the consistency condition,  and thus it follows from the Kolmogorov extension theorem that there exists a probability measure $\tilde\dbP$ on $(\O, \cF_T)$ such that $\{\nu_\pi\}_\pi$ is the finite distribution of the process $Y$ under $\tilde\dbP$. It is straightforward to verify $\tilde\dbP\in \cP(s, \tilde m)$, and $\tilde\dbP_{Y_r} = \dbP_{Y_r}$ for all $r\in [s, T]$. Thus, 
\beaa
J(t, \dbP) =  \int_t^s F(r, \dbP_{Y_r})dr + J(s, \dbP) =  \int_t^s F(r, \dbP_{Y_r})dr +J(s, \tilde \dbP) \le  \int_t^s F(r, \dbP_{Y_r})dr + V(s, {\dbP}_{Y_{s-}}).
\eeaa
 Since $\dbP\in \cP(t,m)$ is arbitrary, we obtain $ V(t,m) \le \tilde V(t,m) $. 

On the other hand,  given $\tilde m$, by Proposition \ref{existence} there exists $\tilde\dbP\in \cP(s, \tilde m)$ such that $J(s, \tilde \dbP) = V(s, \tilde m)$. For the above time partition $\pi$, we introduce another finite measure: for any $A_i\in \cB(\bS)$,
\beaa
\nu_\pi(A_0\times \cds \times A_{m+n}) := \int_{\bS}\dbE^\dbP\Big[  \prod_{i=0}^m \1_{A_i}(Y_{t_i}) \Big|Y_{s-} = y\Big] \times \dbE^{\tilde \dbP}\Big[ \prod_{j=1}^n \1_{A_{m+j}}(Y_{t_{m+j}}) \big| Y_{s-}=y\Big] \tilde m(dy). 
\eeaa
Applying the Kolmogorov extension theorem again there exists a probability measure $\hat\dbP$ on $(\O, \cF_T)$ such that $\{\nu_\pi\}_\pi$ is the finite distribution of the process $Y$ under $\hat\dbP$. It is clear that $\hat\dbP = \dbP$ on $\cF_{s-}$, and $\{Y_{s-}, Y_r, s\le r\le T\}$ has the same distribution under $\hat \dbP$ and $\tilde \dbP$, and $\{Y_r, r<s\}$ and $\{Y_r, r\ge s\}$ are conditionally independent under $\hat \dbP$, conditional on $Y_{s-}$. We shall emphasize that this conditional independence is valid only conditional on $Y_{s-}$, the process $Y$ is in general not Markov under $\hat\dbP$. It is obvious that the processes $M, N$ in \eqref{martingalepb} remain to be $\hat\dbP$-martingales on $[t,s]$. Moreover, for any $s\le s_1 <s_2\le T$, any $0=t_0<\cds< t_m <s \le t_{m+1} <\cds <t_{m+n}=s_1$, and any bounded measurable function $\f_1: \dbR^{(m+1)d}\to \dbR$, $\f_2:\dbR^{nd}\to \dbR$,
\beaa
&&\dis \dbE^{\hat \dbP}\Big[ [M_{s_2}- M_{s_1}] \f_1(Y_{t_0},\cds, Y_{t_m}) \f_2(Y_{t_{m+1}}, \cds, Y_{t_{m+n}})\Big]\\
&&\dis= \dbE^{\hat\dbP}\Big[ \dbE^{\tilde \dbP}\big[[M_{s_2}- M_{s_1}]  \f_2(Y_{t_{m+1}}, \cds, Y_{t_{m+n}}) \big| Y_{s-}\big]\times \dbE^\dbP\big[ \f_1(Y_{t_0},\cds, Y_{t_m}) \big| Y_{s-}\big]\Big]\\
&&\dis = \dbE^{\hat\dbP}\Big[ 0\times \dbE^\dbP\big[ \f_1(Y_{t_0},\cds, Y_{t_m}) \big| Y_{s-}\big]\Big] =0.
\eeaa
Then $M$ is a $\hat \dbP$-martingale on $[s, T]$ as well, and hence a  $\hat \dbP$-martingale on $[t, T]$. Similarly we can show that $N$ is a  $\hat \dbP$-martingale on $[t, T]$, then $\hat \dbP\in \cP(t, m)$. Therefore,
\beaa
\int_t^s F(r, \dbP_{Y_r})dr  + V(s, {\dbP}_{Y_{s-}}) &=& \int_t^s F(r, \dbP_{Y_r})dr  + J(s, \tilde \dbP)\\
&=&\int_t^s F(r, \hat\dbP_{Y_r})dr  + J(s, \hat \dbP) = J(t, \hat \dbP) \le V(t, m).
\eeaa
Since $\dbP\in \cP(t,m)$ is arbitrary, we obtain $ \tilde V(t,m) \le  V(t,m)$, and hence the first equality of \reff{weakDPP}.   
  
It remains to prove the second equality of \reff{weakDPP}. First, since $I_s\le I_{s-}$, it is obvious that $\cP(s, \dbP_{Y_s}) \subset \cP(s, Y_{s-})$, and thus $V(s,  \dbP_{Y_s})\le V(s,  \dbP_{Y_{s-}})$ for all $\dbP\in \cP(t, m)$. On the other hand, for any $\dbP\in \cP(t, m)$, set $\tilde \dbP\in \cP(t, m)$ be such that $\tilde \dbP = \dbP$ on $\cF_{s-}$ and $I_r = I_{s-}$ for all $r\ge s$, $\tilde \dbP$-a.s. Then $\tilde \dbP_{Y_s} = \dbP_{Y_{s-}}$, and thus 
\beaa
 \int_t^s F(r, \dbP_{Y_r})dr + V(s, \dbP_{Y_{s-}}) =  \int_t^s F(r, \tilde\dbP_{Y_r})dr + V(s, \tilde \dbP_{Y_{s}}) \le  \sup_{\dbP \in \cP(t,m)}  \int_t^s F(r, \dbP_{Y_r})dr  + V(s, {\dbP}_{Y_{s}}).
 \eeaa
 This completes the proof immediately.
 \qed

 In order to derive the dynamic programming equation, we  follow the usual procedure, which requires Itô's formula along the flow of measures $\{\dbP_{Y_s}\}_{t \le s \le T}$, as we shall develop in the next section.

\section{Itô's formula for flows of laws of semimartingales}\label{Ito}

In contrast with the available literature reviewed in the introduction,
our Itô's formula allows for possible jumps for both the semimartingale and its flow of marginal laws $\bm = \{m_s\}$. The mapping $s\mapsto m_s$ is also {\cad} and we shall denote
\bea
\label{JAm}
J_\dbT(\bm) := \{s \in \dbT: m_s \neq m_{s-}\},\q J^c_\dbT(\bm) := \{s \in \dbT: m_s = m_{s-}\},\q \mbox{for all} ~\dbT\subset [0, T].
\eea

We first introduce the notion of \textit{linear functional derivative}, in the same spirit as Carmona \& Delarue \cite[Vol 1, Definition 5.43]{CarDel} and Cardialaguet, Delarue, Lasry \&  Lions \cite{CDLL}:
\begin{defn}\label{linderiv}
{\rm (i)} $u : \cP_2(\dbR^{d'}) \longrightarrow \dbR$ has a linear functional  derivative if there exists
$$ \d_m u : \cP_2(\dbR^{d'}) \times \dbR^{d'} \rightarrow \dbR $$
such that  $\d_m u$ is continuous for the product topology and

$\bullet$  the mapping $y \mapsto \d_m u(m,y)$ has quadratic growth in $y$, locally uniformly in $m$. That is,  for any compact set $\Xi\subset \cP_2(\dbR^{d'})$, $\sup_{m\in \Xi} |\d_m u(m, y)|\le C_\Xi [1+|y|^2]$.

$\bullet$ for all $m,m' \in \cP_2(\dbR^{d'})$,
\bea
\label{dmu}
u(m')-u(m) = \displaystyle\int_0^1 \int_{\dbR^{d'}} \d_m u(\l m' + (1-\l)m,y)(m'-m)(dy)d\l.
\eea
{\rm (ii)} $C_2^{1,2}([0, T]\times \cP_2(\dbR^{d'}))$ denotes the set of functions $u: [0, T]\times \cP_2(\dbR^{d'})\to \dbR$ such that 

$\bullet$  $\pa_t u$, $\d_m u$, $\pa_y \d_m u$, $\pa_{yy}^2 \d_m u$ exist and are continuous in all variables; 

$\bullet$ $\pa_{yy}^2 \d_m u$ is bounded in $y$, locally uniformly in $(t, m)$.
\end{defn}
Here the subscript $_2$ in $C^{1,2}_2$ is to refer the growth conditions so as to ensure appropriate square integrability in the analysis below.

By abusing the notation, in the following statement, we let $Y$ denote a general càdlàg $\dbR^{d'}-$valued semimartingale on $[0, T]$. We denote $Y^c$  the continuous part of $Y$;  $Y^c_t = Y_0 + M^c_t + A^c_t$ the Doob-Meyer decomposition, where $M^c$ is the martingale part and $A^c$ is the finite variation part;   $\|A^c\|_t$  the total variation process of $A^c$ and $\la M^c\ra_t$ the quadratic variation process of $M^c$. 
\begin{thm}[Itô's formula]\label{Itothm}
Let $u\in C^{1,2}_2([0, T]\times \cP_2(\dbR^{d'}))$, and assume 
\bea\label{summability}
 \dbE\Big[ \|A^c\|_T^2 + \la M^c\ra_T  + \Big(\sum_{0 < s \le T} \lvert  Y_s-Y_{s-} \rvert \Big)^2 \Big] < \infty.
\eea
Then, denoting $\bm = \{m_s\}_{0\le s\le T}$  the marginal laws of $Y_s$,
\bea
\label{Ito0}
&&u(T,m_T)= u(0,m_0) + \int_0^T \pa_t u(s,m_s)ds \nonumber \\
&&\q+ \dbE\Big[\int_0^T \pa_y \d_m u(s,m_s, Y_s)\cd dA_s^c + \frac{1}{2}\int_0^T \pa_{yy}^2\d_m u(s,m_s,Y_s):d\langle M^c \rangle_s\Big] \\
&&\q +  \!\!\!\!\sum_{s \in J_{(0, T]}(\mathbf{m})} \!\!\!\! [u(s, m_s) -u(s, m_{s-})]  +  \dbE\Big[ \sum_{s\in J_{(0, T]}^c(\bm)}  \!\!\!\!\big(\d_m u(s, m_s, Y_s)- \d_m u(s, m_s, Y_{s-})\big)\Big]. \nonumber
\eea
\end{thm}

 The proof of this result  is relegated to Appendix \ref{appendixB}. Note that \reff{Ito0} exhibits two different sums: one refers to the jumps of $Y$, while the other to the jumps of the marginals $\bm$. The Poisson process provides a simple example of pure jump process with continuous marginals (i.e., $J_{(0, T]}(\mathbf{m}) = \emptyset$). 
 
 \begin{rem}\label{rem:gpw} {\rm
 The above Itô's formula was derived independently by Guo, Pham \& Wei \cite{GPW} by using a density argument through cylindrical functions. Our approach is more straightforward, as it reduces quickly to the proof of the standard Itô's formula. Notice that our set of conditions is slightly different from theirs (none of them implies the other), see Remark 3.14 in \cite{GPW}. Notice also that we may have stated our results under different sets of assumptions, as the proof requires appropriate integrability conditions on the product between the derivatives of $u$ and the corresponding characteristics of the semimartingale $Y$. Clearly, this can be achieved by a trade-off between the conditions on $u$ and $Y$.  
 \qed}
 \end{rem}
 
 We now specialize the discussion to the case $Y := (X,I)$. Note that $\cP_2(\bS) \subset \cP_2(\dbR^{d+1})$, we may restrict Definition \ref{linderiv} to $\cP_2(\bS)$ only.

 \begin{defn}\label{C12S}
Let $C^{1,2}_2(\ol \bQ_0)$ denote the set of functions $u: \ol \bQ_0 \to \dbR$ such that $\pa_t u, \d_m u, \pa_x \d_m u$, $\pa_{xx}^2\d_m u$ exist and are continuous in all variables, and $\pa_{xx}^2\d_m u$ is bounded in $x$,  locally uniformly in $(t, m)$, where the functional linear derivative takes the form $\d_m u: (t, m, x, i) \in \ol \bQ_0\times \dbR^d \times \{0,1\}\to \dbR$ satisfying, for any $t\in [0, T]$ and $m,m' \in \cP_2(\bS)$,
\beaa
u(t,m')-u(t, m) = \displaystyle\int_0^1 \int_{\bS} \d_m u(t, \l m' + (1-\l)m, x, i)(m'-m)(dx, di)d\l.
\eeaa
\end{defn}
In this case, of course there is no need to consider the derivative of $\d_m u$ with respect to the $i$-variable. Instead, we denote 
 \bea
 \label{DIu}
 \d_m u_i(t,m,x) := \d_m u(t,m,x,i) \ \mbox{for $i \in \{0,1\}$, and} \ D_I u := \d_m u_1 - \d_m u_0.
 \eea
 
 \begin{eg}
Let us define, for a given probability measure $\dbP$,
$ u(m) := \f(m[\psi]),$
with $\psi$ smooth and 
$m[\psi] :=  \sum_{i = 0, 1}\int_{\dbR^d}\psi(x,i)m(dx,i).$
Then we compute
$$ \d_m u(m,x,i) = \f'(m[\psi])\psi(x,i) \ \mbox{and} \ D_I u(m,x) = \f'(m[\psi])[\psi(x,1)-\psi(x,0)].$$
 \end{eg}

 Recall the infinitesimal generator of $X$, we define
\bea\label{DiffOp}
\left.\ba{c}
\dis\dbL u(t,m) :=  \pa_t u(t,m) + \int_{\dbR^d}  \cL_x\d_m u_1(t,m,x) m(dx,1),\q \mbox{where}\\
\dis \cL_x\d_m u_1(t,m,x):= b(t,x,m)\cd \pa_x \d_m u_1(t,m,x) + {1\over 2} \si^2(t,x,m): \pa_{xx}^2  \d_m u_1(t,m,x). 
\ea\right. 
\eea
We now state the It\^{o} formula for  $\bm := \{m_s:= \dbP_{Y_s}\}_{s\in [-1, T]}$. Note that in Theorem \ref{Itothm}, we consider the jumps on $(0, T]$. However, in light of DPP \reff{weakDPP}, it is more convenient to consider the jumps on $[0, T)$, namely we  include the jump at the initial point instead of the ending point. Such an adjustment is straightforward.
 
 \begin{cor}\label{CorIto}
Let $m\in \cP_2(\bS)$, $\dbP\in \cP(0, m)$, and $u \in C^{1,2}_2(\ol \bQ_0)$.  Then, 
 \bea
 \label{JumpFunIto}
 \left.\ba{lll}
\dis  u(T,m_{T-}) = u(0,m) + \int_0^T \dbL u(s,m_s)ds \\ 
\dis\qq  + \sum_{s \in J_{[0, T)}(\mathbf{m})} [u(s,m_s)-u(s, m_{s-})] + \dbE^\dbP\Big[\int_{J^c_{[0, T)}(\mathbf{m})} D_I u(s,m_s,X_s)dI_s \Big].
\ea\right.
\eea
\end{cor}
\proof We can easily see that $Y_s-Y_{s-} = (0, I_s-I_{s-})$ and 
\beaa
Y^c_s = (X_s, I_{0-}),\q dM^c_s=  (\si(s, X_s, m_s) dW^\dbP_s, 0),\q dA^c_s = (b(s, X_s, m_s) ds, 0).
\eeaa
Then \reff{summability} obviously holds true. 
Now following Theorem \ref{Itothm}, but by considering the jump at $0$ instead of at $T$, we have
\beaa
&&u(T,m_{T-})- u(0,m) = \int_0^T \pa_t u(s,m_s)ds + \dbE^\dbP\Big[\int_0^T \pa_y \d_m u_1(s,m_s, X_s)\cd b(s, X_s, m_s)ds \\
&&\qq + \frac{1}{2}\int_0^T \pa_{xx}^2\d_m u_1(s,m_s,X_s):d\langle X \rangle_s\Big] + \sum_{s \in J_{[0, T)}(\mathbf{m})}\!\!\! \!\!\!\! [u(s, m_s) -u(s, m_{s-})]  \\
&&\qq+  \dbE^\dbP\Big[ \!\!\!\! \sum_{s\in J_{[0, T)}^c(\bm)} \!\!\! \!\!\!\!\big[\d_m u(s, m_s, Y_s)- \d_m u(s, m_s,Y_{s-})\big]\Big] \\
 &&= \int_0^T \dbL u(s,m_s)ds +\!\!\!\! \sum_{s \in J_{[0, T)}(\mathbf{m})}\!\!\!\! [u(s,m_s)-u(s, m_{s-})]+ \dbE^\dbP\Big[\int_{J^c_{[0, T)}(\mathbf{m})} \!\!\!\!\!\!\!\!D_I u(s,m_s,X_s)dI_s \Big],
 \eeaa
 where the last equality thanks to the fact that $I_s \neq I_{s-}$ if and only if $I_s = 0, I_{s-}=1$.
\qed

We remark that, in this case $J_{[0, T]}(\mathbf{m}) = \{s \in [0,T] : \dbP(\t = s) > 0\}$. That is,  $J_{[0, T]}(\mathbf{m})$ is the collection of all atoms of $\t$ under $\dbP$.

\section{Obstacle problem on the Wasserstein space}\label{obs}

\subsection{The dynamic programming equation}

We first introduce a partial order $\preceq$ on $\cP_2(\mathbf{S})$: we say that $m' \preceq m$ if
\bea\label{order}
m'(dx, 1) = p(x) m(dx, 1), \ \mbox{and} \ m'(dx,0) = [1-p(x)]m(dx,1) + m(dx, 0), 
\eea
for some measurable $p : \dbR^d \rightarrow [0,1]$, i.e. $m'(dx,1)$ is obtained from $m$ by randomly stopping a proportion $1-p(x)$ of the surviving particles. In our context, $m_{t^-} = \dbP_{Y_{t-}}$ and $m_t = \dbP_{Y_t}$, with $\dbP \in \cP(t,m)$, so that $m_t \preceq m_{t^-}$ with conditional transition probability 
\bea
\label{ptx}
p(x) = p(t,x) := \dbP(I_t = 1 \mid X_t = x, I_{t-} = 1).
\eea

\begin{rem}\label{compact} {\rm
The set $\{m' : m' \preceq m\}$ is compact, as it is in continuous bijection with $\{ \hat m \in \cP_2(\mathbf{S} \times \{0,1\}) : \hat m \circ (\mathbf{x}, \mathbf{i})^{-1} = m \}$, with $(\mathbf{x}, \mathbf{i}, \mathbf{i'})$ the projection coordinates on $\mathbf{S} \times \{0,1\}$. }
\end{rem}
Our main objective is to show that the dynamic programming equation corresponding to our mean field optimal stopping problem, as deduced from the DPP \eqref{weakDPP}, is
\bea\label{obstacle}
\left.\ba{c}
\dis \underset{m' \in C_u(t,m)}{\min}[-(\dbL u + F)(t,m')] = 0, 
~D_I u(t,m,\cd) \ge 0,~ u(T,\cd) =g, ~\mbox{for all } (t,m) \in \bQ_0,\\
\dis \mbox{where}\q C_u(t,m) := \Big\{m' \preceq m : u(t,m') = u(t,m)\Big\}.
\ea\right.
\eea
  By analogy with standard optimal stopping, we call \eqref{obstacle} \textit{obstacle problem} on the Wasserstein space. 
The different components of this equation have the following interpretation. \qed

\begin{rem}
\label{rem-Obstacle}
\no{\rm (i)  As will be proved in Lemma \ref{monotonic}, the inequality $D_I u(t,m,\cd) \ge 0$
 expresses the natural monotonicity of the optimal stopping problem, i.e. $u$ is increasing for $\preceq$. In other words, the larger the set of surviving particles is, the larger the value function is. \\ 
(ii) $C_u(t,m)$ is the collection of \textit{admissible stopping strategies at time $t$}, i.e. those that preserve the value function for smaller sets of surviving particles. 
\\
(iii) The equation
 $ \underset{m' \in C_u(t,m)}{\min}\!\!\!\!-(\dbL u + F)(t,m') = 0 $
 characterizes the sets of particles that are optimal to keep diffusing (in the same spirit as the classical HJB equation, where the min characterizes the optimal controls). Note that $C_u(t,m)$ is compact, as a closed subset of the compact set $\{m' \preceq m \}$, see Remark \ref{compact}. Therefore the $\min$ is attained by the continuity of $(\dbL u + F)(t, \cd)$. Finally, as $m \in C_u(t,m)$, we have
$-(\dbL u + F)(t,m) \ge 0$. \\
(iv) The boundary condition $ u(T,\cd) =g$ is due to \reff{weakoptstop} directly. Moreover, the boundary condition implies that $u(t, m) = g(m)$, for all $t\in [0, T]$ and $m\in \pa \cP_2(\bS) := \{ m \in \cP_2(\bS) : m(\dbR^d, 1) = 0 \}$, i.e. all particles are stopped. Indeed, in this case $\{m': m' \preceq m\} = \{m\}$ and thus $C_u(t, m) = \{m\}$. Recall \reff{F} and \reff{DiffOp}, then \eqref{obstacle} implies to $-\pa_t u(t,m)= -(\dbL u+F)(t,m) = 0$. This clearly implies that $u(t,m) = u(T,m) = g(m)$ for all $t\in [0, T]$. 
\qed}
  \end{rem}

\begin{lem}\label{monotonic}
Let $u : \cP_2(\mathbf{S}) \rightarrow \dbR$  admit a linear functional derivative. Then $u$ is nondecreasing for $\preceq$ if and only if $D_I u(m,\cd)\ge 0$ for all $m \in \cP_2(\mathbf{S})$.
\end{lem}

\proof  First, assume  $D_I u(m,\cd)\ge 0$ for all $m \in \cP_2(\dbR^d)$. Then, for $m' \preceq m$ with corresponding transition probability $p$, we have
\beaa
u(m)-u(m') =  \int_0^1\int_{\dbR^d} D_I u(\l m + (1-\l) m', x)[1-p(x)] m(dx,1)d\l \ge 0.
\eeaa

Conversely,  assume that $u$ is nondecreasing for $\preceq$, i.e. $u(m') \preceq u(m)$ for all $m' \preceq m$. Introduce $\cN := \{ x : D_I u(m,x) < 0 \}$; $p_\e(x) := 1-\e \1_\cN(x)$, $\e \in (0,1)$; and the corresponding measure $m'_\e$ defined by \eqref{order}. Then  $(m-m'_\e)(dx,di) =  (2i-1)\e\1_\cN(x)m(dx,1)$, and thus
\bea
\label{cNintegral}
0 &\le& {1\over \e}[u(m)-u(m')] = {1\over \e}\int_0^1\int_{\mathbf{S}} \d_m u\big( \l m +[1- \l]m'_\e, x, i\big)(m-m'_\e)(dx,di)d\l \nonumber\\
&=& \int_0^1 \int_{\cN}D_I u( \l m +[1- \l]m'_\e, x) m(dx,1)d\l.
\eea
Note that $\{\l m +[1- \l]m'_\e: \l\in [0, 1], \e\in [0,1]\}\subset \cP_2(\bS)$ is compact, then $D_I u( \l m +[1- \l]m'_\e, x)$  has quadratic growth in $x$, uniformly in $\l,\e$. Moreover,  sending $\e\to 0$, since $m'_\e\to m$ and $D_I u$ is continuous in $m$,  applying the dominated convergence theorem we obtain from \reff{cNintegral} that $\int_{\cN}D_I u(m, x) m(dx,1) \ge 0$, which is possible only if $m(\cN,1) = 0$. That is, $D_I u(m, x) \ge 0$ for $m(\cd, 1)$-a.e. $x$. Since $D_I u$ is continuous in $(m, x)$ and the set $\{m\in \cP_2(\bS): \supp(m(\cd, 1))=\dbR^d\}$ is dense in $\cP_2(\bS)$, then one can easily show that $D_I u(m, x) \ge 0$ for all $(m,x)\in \cP_2(\bS)\times \dbR^d$.
\qed

\subsection{The main results}

\begin{thm}\label{NecCond}
If $V$ defined in \eqref{weakoptstop} is in $C^{1,2}_2(\ol \bQ_0)$,  then it is a solution of \eqref{obstacle}. 

Moreover, for any $(t,m) \in [0,T) \times \cP_2(\mathbf{S})$, $\dbP^* \in \cP(t,m)$ is optimal for $V(t,m)$  if and only if,  denoting  $\bm^*=\{m_s^* := \dbP^*_{Y_s}\}_{0\le s\le T}$, 
\bea\label{optimality}
\left.\ba{c}
\dis -(\dbL V + F)(s,m_s^*) = 0,\q V(s, m^*_s) = V(s, m^*_{s-}),\q \mbox{for all $s \in [t,T]$}, \\ 
 \dis D_I V(\t,m_\t^*,X_\t)\1_{J^c_{[t,T)}(\bm^*)}(\t) = 0, \q \dbP^*\!\!\mbox{-a.s.,  where}~ m^*_\t := m^*_s|_{s=\t}. 
\ea\right.
\eea
\end{thm}
\proof  \textit{Step 1:} We first prove that
\bea\label{supersolution}
-(\dbL V + F)(t,m) \ge 0, \q D_I V(t,m,x) \ge 0,\q  \mbox{for all $(t,m, x) \in \bQ_0 \times \dbR^d$.}
\eea
Fix $(t, m)$. For any $m' \preceq m$, we may choose $\dbP \in \cP(t,m)$ such that $\dbP_{Y_t} = m'$, and $I_s = I_t$, $s \in [t,T)$, $\dbP$-a.s.  
 For $\d \in (0,T-t)$, by DPP \eqref{weakDPP} we have 
 $$V(t,m) \ge  \int_t^{t+\d} F(s, \dbP_{Y_s})ds  + V(t+\d, \dbP_{(X_{t+\d}, I_t)}).$$
Send $\d\to 0$, note that $\dbP_{(X_{t+\d}, I_t)} \to \dbP_{(X_t, I_t)} = m'$. Then by the continuity of $V$ we have $V(t,m) \ge V(t, m')$.
Since $m' \preceq m$ is arbitrary, by Lemma \ref{monotonic} we see that $D_I V \ge 0$.

To prove that $-(\dbL V +F)(t,m) \ge 0$, we consider $\dbP\in \cP(t, m)$ such that $I_s = I_{t-}$, $s\in [t, T]$, $\dbP$-a.s. Apply Itô's formula \eqref{JumpFunIto} on $[t, t+\d]$ under $\dbP$, we see that all the terms involving the jumps are equal to $0$. Then by DPP \reff{weakDPP} we have, denoting $m_s:=  \dbP_{Y_s}$,
\beaa
0 \ge V(t+\d, m_{(t+\d)-}) - V(t, m) + \int_t^{t+\d} F(s, m_s)ds = \int_t^{t+\d} (\dbL u + F)(s, m_s) ds.
\eeaa
Note that $m_s\to m$ as $s\downarrow t$. Then by the continuity of $\dbL u + F$ one can easily see that $-(\dbL V + F)(t,m) \ge 0$.

\textit{Step 2:} In this step we prove the equivalence of the optimality condition \reff{optimality}. First, if $\dbP^*\in \cP(t, m)$ satisfies \reff{optimality}, applying Itô's formula \eqref{JumpFunIto} on $[t, T)$ we obtain immediately 
\beaa
V(t,m) = V(T, \dbP^*_{Y_{T-}}) - \int_t^T \dbL V(s, \dbP^*_{Y_s})ds  =  g(\dbP^*_{Y_T}) +\int_t^T F(s, \dbP^*_{Y_s})ds.
\eeaa
 As $f$ has quadratic growth in $x$, locally uniformly in $(t,m)$, we may switch the integral and the expectation in the expression of $F$, and thus $\dbP^*$ is optimal.

On the other hand, for any optimal $\dbP^* \in \cP(t,m)$ such that
$ V(t,m) =  \int_t^T F(s, \dbP^*_{Y_s})ds  + g(\dbP^*_{Y_T}).$ 
Denoting $m_s^* := \dbP^*_{Y_s}$, $s\ge t$, with $m = m_{t-}^*$, then by DPP \eqref{weakDPP} and Itô's formula \eqref{JumpFunIto} we have 
\beaa
0&=&\hspace{-2mm} \int_t^T(\dbL V + F)(s,m_s^*)ds  +\hspace{-4mm} \sum_{s \in J_{[t,T)}(\mathbf{m^*})}\hspace{-3mm}[V(s,m_s^*) - V(s, m^*_{s-})] + \dbE\Big[\int_{J^c_{[t,T)}(\mathbf{m^*})}\hspace{-8mm} D_IV(s,m_s^*,X_s)dI_s\Big].
\eeaa
By  Step 1 we have $V(s,m_s^*) \le V(s, m^*_{s-})$. Together with \eqref{supersolution},  we see that all the three terms in the right side above are nonpositive, then all of them should be $0$:
\bea
\label{optimality2}
\left.\ba{c}
\dis (\dbL V + F)(s,m_s^*) =0,~ \mbox{a.e.}~s\in [t, T];\q V(s,m_s^*) = V(s, m^*_{s-}),~ \mbox{for all } s\in J_{[t,T)}(\mathbf{m^*});\\
\dis  \int_{J^c_{[t,T)}(\mathbf{m^*})}\hspace{-8mm} D_IV(s,m_s^*,X_s)dI_s=0,~\dbP^*-\mbox{a.s.}
\ea\right.
\eea
Since $\dbL V + F$ is continuous and, for $s\in J^c_{[t,T)}(\mathbf{m^*})$,  by definition $m^*_s = m^*_{s-}$ and hence $V(s,m_s^*) = V(s, m^*_{s-})$, then the first line of \reff{optimality2} implies that the first line of \reff{optimality} holds for all $s\in [t, T)$. Moreover, since $\t$ is the only jump point of $I$, the second line of  \reff{optimality2} is clearly equivalent to the second line of \reff{optimality}.

\textit{Step 3:} Finally we complete the verification of \reff{obstacle}. First by \reff{weakoptstop} $V(T,m)=g(m)$. Then, by Step 1, it remains to verify $\min_{m'\in C_V(t,m)}[ - (\dbL V + F)(t, m')]=0$. Note that $m^*_{t-}=m$ and set $s=t$ in the first line of \reff{optimality}, we have $m^*_t \in C_V(t, m)$. Thus \beaa
0\le\min_{m'\in C_V(t,m)}[ -(\dbL V + F)(t, m')] \le -(\dbL V + F)(t, m^*_t) =0,
\eeaa
 and therefore the equality holds.
\qed

\begin{thm}[Verification]\label{verification}
Let $u\in C^{1,2}_2(\ol \bQ_0)$ be a solution of \eqref{obstacle}. Then $u = V$.  
\end{thm}
\proof We prove the theorem by using the obstacle equation \eqref{obstacle} to construct an $\e$-optimal control for \eqref{weakoptstop}. We fix $m\in \cP_2(\bS)$ and assume for simplicity that $t=0$.

\textit{Step 1:} We first prove that $u \ge V$. For an arbitrary $\dbP \in \cP(0,m)$, we apply Itô's formula \reff{JumpFunIto}  and obtain: again denoting $\bm = \{m_s:= \dbP_{Y_s}\}$,
\bea
 \label{verif1}
 \left.\ba{lll}
\dis  u(T,m_{T-}) = u(0,m) + \int_0^T \dbL u(s,m_s)ds \\ 
\dis\qq  +\sum_{s \in J_{[0, T)}(\mathbf{m})} [u(s,m_s)-u(s, m_{s-})] + \dbE^\dbP\Big[\int_{J^c_{[0, T)}(\mathbf{m})} \!\!\!\! \!\!\!\!D_I u(s,m_s,X_s)dI_s \Big].
\ea\right.
\eea
By $\eqref{obstacle}$ and Lemma \ref{monotonic} we have $u(s, m_s) \le u(s, m_{s-})$. Then,
\reff{obstacle} and \eqref{verif1} imply that 
\beaa
u(0,m) \ge u(T,m_{T-}) - \int_t^T \dbL u(s,m_s)ds \ge  g(m_T)+ \int_t^T F(s,m_s)ds.
\eeaa
 Since $\dbP\in \cP(0, m)$ is arbitrary, we obtain  $u(0,m) \ge  V(0,m)$.

\textit{Step 2:} We now show that $u \le V$. Let $n\ge 1$, $t_j := {j\over n}T$, $j=0,\cds, n$. We define $\dbP^n\in \cP(0, m)$ and $m_s^n := \dbP^n_{Y_s}$ recursively such that $m^n_{0-} = m$, and for $j=0,\cds, n-1$,  thanks to Remark \ref{compact},
$$ m^n_{t_j}\in C_u(t_j, m_{t_j-}^n) \ \mbox{s.t.} \ -(\dbL u + F)(t_j, m^n_{t_j}) = 0 \q \mbox{and}\q  m_s^n \circ \mathbf{i}^{-1} = m^n_{t_j} \circ \mathbf{i}^{-1}, ~s \in [t_j,t_{j+1}).
$$
By the arguments of Proposition \ref{existence} applied to $\cP(0, m)$, and by \reff{asympt}, one can easily show that $\cW_2(m_s^n, m_{t_j}^n) \le {C_{m}\over \sqrt{n}}$, $s \in [t_j, t_{j+1})$, for some constant $C_{m}>0$ which may  depend on  $m$ but is uniform on $n$. Moreover, by Proposition \ref{existence} and the compactness of $[0, T]$, we see that the set $\{\dbP_{Y_s}, \dbP_{Y_{s-}}: s\in [0, T], \dbP\in \cP(0, m)\}$ is compact. 
As $u \in C^{1,2}_2(\ol \bQ_0)$, $\dbL u+F$  is continuous and then uniformly continuous on this set. Then, there exists a modulus of continuity function $\rho$ such that 
\beaa
-(\dbL u + F)(s, m_s^n)  = -(\dbL u + F)(s, m_s^n) + (\dbL u + F)(t_j, m^n_{t_j}) \le \rho({T\over n} + {C_m\over \sqrt{n}}),\q s\in [t_j, t_{j+1}).
\eeaa
By Itô's formula \reff{JumpFunIto}, and noting that  $\dbP^n$ is constructed such that there is no contribution of the jump terms, we have
\beaa
u(0,m) &=&  u(T,m_{T-}^n)  - \int_0^T \dbL u(s, m^n_s) ds \\
&\le& g(m^n_T) + \int_t^T F(s,m_s^n)ds + T \rho({T\over n} + {C_m\over \sqrt{n}}) \le V(0, m) +  T \rho({T\over n} + {C_m\over \sqrt{n}}).
\eeaa
Send $n\to \infty$, we obtain $u(0, m)\le V(0, m)$.
\qed

\subsection{Some discussions on optimal stopping policies} 

Proposition \ref{existence} (i) guarantees that the mean field optimal stopping problem has an optimal randomized stopping strategy, i.e. a probability measure $\dbP^*$ on $\O$ s.t. $g(\dbP^*_{Y_T}) +\int_t^T F(s, \dbP^*_{Y_s}) ds= V(t,m)$.
  A pure stopping strategy corresponds to the case where the conditional transition probability in \reff{ptx} $p_s(\cd) \in \{0,1\}$ for all $s \in [t,T]$. In this case, the optimal stopping time is in closed-loop, i.e. $\t$ is a stopping time w.r.t to the $\dbP^*$-augmented filtration of $X$, and the obstacle equation \eqref{obstacle} reduces to:
\bea\label{obstacleX}
\underset{A \in \cB_u(t,m)}{\min}[-(\dbL u+F)(t,m^A)] \!=\! 0, 
~u(t,m)\!=\!\max_{A \in \cB(\dbR^d)}u(t,m^A),~ u|_{t=T} \!=\! g,
~(t,m) \in \ol \bQ_0, ~~~
\eea
where
$m^A := m \circ (\mathbf{x}, \mathbf{i} \1_A(\mathbf{\bx}))^{-1},$ and $\cB_u(t,m) := \big\{A \in \cB(\dbR^d) :  u(t,m^A) = u(t,m)  \big\}.$

We now discuss heuristically how to use the value function $V$ to construct an optimal stopping time, provided $V\in  C^{1,2}_2(\ol \bQ_0)$. In light of \reff{optimality} and recalling that $D_I V\ge 0$, introduce 
\bea
\label{O}
K(t, m) := \big\{x\in \dbR^d: D_I V(t,m, x) =0\big\}.
\eea
Fix $(0, m_{0-})$. We set $m^*_{0-}:= m_{0-}$ and construct $\bm^*$ for $V(0, m_{0-})$ in several steps.

{\it Step 1.} First, by \reff{obstacle} and Remark \reff{compact}, there exists $m^*_0\in C_V(0, m^*_{0-})$ such that $m^*_0 \preceq m^*_{0-}$ and $(\dbL V + F)(0, m^*_{0})= 0$.  In particular, if $m^*_{0-}\circ \bx^{-1}$ is continuous on $\{I_{0-}=1\}$, there exists $A\in \cB(\dbR^d)$ such that $I_0 = I_{0-} \1_{A^c}(X_0)$, and thus the optimal stopping time is a pure strategy at $0$. 

{\it Step 2.} Let $\dbP^*$ be  a weak solution to the following McKean-Vlasov SDE:   
\bea
\label{optimalMVSDE}
\dbP^*_{Y_0} = m^*_0,\q \mbox{$X$ satisfies \reff{WP} and}\q I_t = I_0 \1_{\{X_s \in K(s, \dbP^*_{Y_s}), 0< s\le t\} },\q \dbP^*\mbox{-a.s.}
\eea
Assume $m^*_t:= \dbP^*_{Y_t}$ is continuous up to certain $t_1>0$. Then the optimal stopping time between $[0, t_1)$ is a pure strategy: $\t = \inf\{t\ge 0: X_t \notin K(t, m^*_t)\}$. Note that, since $V$ is the value function, we should have $(\dbL V + F)(t, m^*_{t})= 0$ for $t\in [0, t_1)$. We shall remark though, the McKean-Vlasov SDE \reff{optimalMVSDE} is path dependent and has discontinuous coefficients, so in general it is  hard to solve. Moreover, the case that $t_1=0$ is even more difficult to solve.

{\it Step 3.} We have obtained $m^*_{t_1-}$ from Step 2. As in Step 1, we may find $m^*_{t_1} \in C_V(t_1, m^*_{t_1-})$ at $t_1$ such that $m^*_{t_1} \preceq m^*_{t_1-}$ and $(\dbL V + F)(t_1, m^*_{t_1})= 0$. Then following Step 2 again we can hopefully extend $m^*$ to certain $t_2>t_1$. Repeat the procedure, we may construct $m^*$ on $[0, T]$. 

We emphasize again that this procedure is just to illustrate the idea, in particular, it could be helpful for constructing approximate optimal stopping times, as we saw in Theorem \ref{verification} Step 2. In general it is hard to realize this procedure, in fact, even the existence of classical solution is a very challenging task. Nevertheless, in Subsection \ref{sect-classical} below we will present an example where $V$ is smooth and we can construct the $\tau$ explicitly. We also remark again that, the optimal stopping time in the continuous region constructed in Step 2 above is always a pure strategy, while in the jump region in Step 1 the optimal stopping time could be indeed mixed, but will also be a pure strategy when the distribution of the survival particles at that time is continuous.

\section{Examples}\label{examples}

\subsection{Connection with standard optimal stopping}

In this subsection we consider the case  that $b$ and $\si$ do not depend on the $\cP_2(\mathbf{S})$-valued variable. For a measurable function $\f$, we define the optimal stopping problem
\bea\label{MFstandard}
V(t,m) := \sup_{\dbP \in \cP(t,m)} \dbE^\dbP[\f(X_T)],
&(t,m) \in \ol \bQ_0, &
\eea
We  also introduce  $v(t,x):= V(t, \d_{(x,1)})$ which is related to  the standard obstacle problem
\bea\label{varineq}
\min\{-(\pa_t v +\cL v), v-\f\} = 0, \q v(T,\cd) = \f,\q\mbox{where}\q \cL v:=  b\cd \pa_x v + {1\over 2} \si^2: \pa_{xx}^2 v.
\eea
To be consistent with Definition \ref{C12S}, let $C^{1,2}_2([0, T]\times \dbR^d)$ denote the set of $v\in C^{1,2}([0, T]\times \dbR^d)$ such that $ \pa^2_{xx} v$ is bounded.  This condition can be relaxed in this case though.

\begin{prop}\label{standardprop}
Assume $v\in C_2^{1,2}([0, T]\times \dbR^d)$. Then: 
\\
{\rm (i)} $V(t,m) =  \int_{\mathbf{S}}\big(v(t,x)i + \f(x)(1-i)\big)m(dx,di)$, and $V$ is a classical solution of the corresponding obstacle equation on the Wasserstein space;
\\
{\rm (ii)} The probability measure $\dbP^*$ s.t. $\t = \inf\{s \ge t : v(s,X_s) =\f(X_s) \}$ on $\{I_{t-}=1\}$, $\dbP^*$-a.s., is optimal for the problem $V(t,m)$.  In particular, we see that $\t$ is a pure stopping strategy under $\dbP^*$.
\end{prop}
\proof Denote by $u$ the right-hand side of the expression in (i).
Then $u \in C_2^{1,2}(\ol \bQ_0)$ with 
\beaa
&\pa_t u(t,m) = \int_{\dbR^d}\pa_t v(t,x)m(dx,1), \ \d_m u(t,m,x,i) = v(t,x)i+\f(x)(1-i),\\
&\pa_x \d_m u_1(t,m,x) = \pa_x v(t,x), \ \pa_{xx}^2 \d_m u_1(t,m,x) = \pa_{xx}^2 v(t,x) \q \mbox{for all $(x,i) \in \mathbf{S}.$} 
\eeaa
We then show that $u$ is a solution of the equation \eqref{obstacle}.
First, by \eqref{varineq},
$$ D_I u = v - \f \ge 0 \q \mbox{and} \q  -\dbL u(t,m) = \int_{\dbR^d}- \cL v(t,x)m(dx,1) \ge 0.
$$
 Defining $A_t := \{x : v(t,x) - \f(x) > 0 \}$ and $m^{A_t} := m \circ (\mathbf{x}, \mathbf{i} \1_{A_t}(\mathbf{x}))^{-1}$, we have
$$ u(t,m) - u(t,m^{A_t}) = \int_{A_t^c}[v(t,x)-\f(x)]m(dx,1) = 0, $$
 and therefore $m^{A_t} \in C_u(t,m)$. As $- \cL v (t,x) = 0$, $x \in A_t$, we have
 $-\dbL u(t,m^{A_t}) = 0. $
 Thus, $u$ is a solution of \eqref{obstacle}, and we deduce that $u=V$ by Theorem \ref{verification}. 

To see that (ii) holds, notice that the flow $m^*_s := \dbP^*_{Y_s}$ is s.t. $m_s^* = (m_{s-}^*)^{A_s}$ for all $s \in [t,T]$. Then $\dbP^*$ clearly satisfies \eqref{optimality}, and thus is optimal for $V(t,m)$.
\qed

\subsection{Convex functions of the expectation}\label{cvx-example}

Let $d=1$, $\psi, h,\f : \dbR \rightarrow \dbR$, with $\f$ convex. We consider the optimal stopping problem:
\bea\label{convex}
V(t, m) := \sup_{\dbP \in \cP(t,m)}\Big[ \dbE^\dbP[\psi(X_T)] + \f\big(\dbE^\dbP[h(X_T)]\big)\Big].
\eea
This is an extension of the mean-variance optimal stopping problem. Introducing the convex dual $\f^*(\a) := \sup_{\b \in \dbR}\{\a \b - \f(\b)\}$, we may write 
\bea
V(t,m)  = \sup_{\a \in \dbR} \big[-\f^*(\a) + V_\a(t,m)\big], \ \mbox{with} \ V_\a(t,m) := \sup_{\dbP \in \cP(t,m)} \dbE^\dbP[\psi_\a(X_T) ], \ \psi_\a := \psi+\a h.  \nonumber
\eea
Assuming $u_\a(t,x) := V_\a(t,\d_{(x,1)}) \in C^{1,2}_2([0, T]\times \dbR)$, it follows from 
Proposition \ref{standardprop} that
\bea
\label{Va}
\left.\ba{c}
\dis V_\a(t,m) = \int_\mathbf{S} \big[u_\a(t,x)i + f_\a(x)[1-i]\big]m(dx,di),\\
\dis - \dbL V_\a(t, m^{A_t}) = 0,\q \mbox{where}~ A_t := \{x: u_\a(t,x) >  \f(t,x)\};\\
\dis D_I V_\a(t, m, x) = u_\a(t,x) - \f(t,x) =0, \q \mbox{when}~ u_\a(t,x) =  \f(t,x).
\ea\right.
\eea
Since $\a$ is one dimensional, it is not hard to find $\a^*(t, m)$ s.t. $V(t, m) =V_{\a^*(t, m)}(t, m) -\f^*(\a^*(t,m))$.

Moreover, fix $(t, m)$ and let $\dbP^*\in \cP(t, m)$ be the optimal measure for the problem $V_{\a^*(t, m)}(t, m)$, as constructed in the previous subsection. Then it is obvious that $\dbP^*$ is optimal for $V(t, m)$ as well, and by Proposition \ref{standardprop} (ii), $\t$ is an optimal stopping strategy under $\dbP^*$.

\begin{rem}{\rm
Let $d=1$. Another natural example is the optimal stopping of the expected shortfall:
\bea\label{ESoptstop}
V(t,m) := \inf_{\dbP \in \cP(t,m)} \mathrm{ES}_\a^\dbP(X_T) \q \mbox{for all}
~~(t,m)\in\overline{\mathbf{Q}}_0, \nonumber
\eea
for some fixed $\a \in (0,1)$, where $ \mathrm{ES}_{\a}^\dbP$ denotes the expected shortfall under $\dbP$, i.e., for any r.v. $Z$ with law $\mu$,
\bea\label{eq-RU}
&g(\mu):=\dis\mathrm{ES}_\a^\dbP(Z) := \frac 1 \a \int_0^\a q_\g(Z) d\g  =  \inf_{\b \in \dbR} \Big\{ \b + \frac{1}{1-\a}\int_\dbR (x-\b)^+\mu(dx)\Big\}, \\
&\dis\mbox{where}\q q_\g (Z) := \inf\{z : \mu(Z \le z) > \g \}. \nonumber
\eea
Here the second equality has been established by Rockafellar \& Uryasev \cite{RU}. 
%This implies
%\beaa
%V(t,m) = \inf_{\b\in \dbR} \Big\{ \b + \frac{1}{1-\a}V_\b(t,m)\Big\},\q \mbox{where}\q V_\b(t,m) :=  \inf_{\dbP\in \cP(t,m)} \dbE^\dbP[(X_T-\b)^+].
%\eeaa
It is not clear wether this value function is smooth, so that the result of the current paper does not apply. A similar comment applies to the optimal stopping under probability distortion. These two examples are discussed in our accompanying paper \cite{TTZ2}, which adresses the possible non-smoothness by introducing an appropriate notion of viscosity solution.
\qed}
\end{rem}

\subsection{Construction of a smooth solution}
\label{sect-classical}
In this subsection we construct an example where the obstacle problem indeed has a classical solution. 
First, set $b=0$, $\si=1$, and thus, for any $\dbP \in \cP(t, m)$,
\bea
\label{eg-X}
X_s  = X_t  + W^\dbP_{\t\wedge s} - W^\dbP_t\q\dbP\mbox{-a.s. on $\{I_{t-}=1\}$}.
\eea
Next, let  $a \in C^1([0, T])$ and $\f \in C_2^{1,2}\big([0, T]\times \dbR_+\big)$ be positive functions such that 
\bea
\label{eg-phi}
\pa_x \f(t,x) = 0~\mbox{for}~ x\ge a_t,\q \mbox{and}\q \pa_x \f (t,x)>0~\mbox{for}~  x< a_t.
\eea
 One such example can be $\f(t, x) := e^{-[(a_t-x)^+]^3}$. Moreover, we introduce another positive function $\psi\in C^2( \dbR)$ with bounded derivatives, and set
\bea
\label{eg-u0}
u_0(t, m) := [T-t]\f(t, v_0(m)),\q\mbox{where}\q v_0(m) := \int_\dbR \psi(x) m(dx, 1).
\eea

\begin{prop}
\label{prop-eg}
Under the above setting, $u_0 \in C_2^{1,2}([0, T]\times \cP_2(\bf S))$, and $u_0$ is the classical solution to the obstacle problem \reff{obstacle} with
\bea
\label{eg-f}
F(t,m) := - \dbL u_0(t, m) -  \big[v_0(m) -  a_t\big]^+,\q g:= 0.
\eea
\end{prop}
We remark that this $F$ may not take the specific form of \reff{F}, which is mainly motivated from applications but not really required for our theory. Since this example is just for illustration purpose of the theory, we content ourselves by allowing for this more general $F$. We emphasize again that in general it is hard to have classical solution for our obstacle problem, and therefore we shall investigate viscosity solutions in our accompanying paper \cite{TTZ2}. 

\proof  First, by Definition \ref{linderiv} one may easily verify: $\d_m v_0(m, x, 1) = \psi(x)$, $\d_m v_0(m, x, 0) = 0$, $\pa_t u_0(t,m) = [T-t]\pa_t  \f(t, v_0(m)) -  \f(t, v_0(m))$, $\d_m u_0(t, m, x, 1) = [T-t]\pa_x   \f(t, v_0(m))\psi(x)$, $\d_m u_0(t, m, x, 0)=0$.
Then it is clear that $u_0 \in C^{1,2}_2(\ol \bQ_0)$.

We now show that $u_0$ satisfies \reff{obstacle}. Clearly, $u_0(T, .) = 0 = g$, and
\bea
\label{eg-DIu0}
\dis D_I u_0(t, m, x) =  [T-t]\pa_x   \f(t, v_0(m))\psi(x) \ge 0,\q  -  \dbL u_0(t, m) - F(t, m) =  \big[v_0(m) -  a_t\big]^+ \ge 0.
\eea
In particular, $- (\dbL u_0+F)(t, m)=0$ when $v_0(m) \le a_t$.  Finally,  when $v_0(m) > a_t$, combining \reff{eg-phi} and \reff{eg-u0}, we have \bea
\label{eg-Cu0}
C_{u_0}(t, m) =  \{ m' \preceq m: v_0(m') \in [a_t, v_0(m)]\},\q t<T.
\eea
 Set $m'_* \preceq m$ by \reff{order} with $p(x) \equiv {a_t \over v_0(m)}$. Then $m'_*\in C_{u_0}(t, m)$ with $v_0(m'_*) = a_t$.  Therefore,
\beaa
\underset{m' \in C_{u_0}(t,m)}{\min}-(\dbL u_0 + F)(t,m') \le -(\dbL u_0 +F)(t,m'_*)= \big[v_0(m'_*) -  a_t\big]^+ =0.
\eeaa
This, together with \reff{eg-DIu0}, completes the proof.
\qed

In the rest of this subsection, we construct an optimal $\dbP^*\in \cP(0, m_{0-})$ for the problem $V_0:= V(0, m_{0-}) = u_0(0, m_{0-})$. For simplicity we assume
\bea
\label{eg-psi}
T=2,~ \psi(x) := e^{-{x^2\over 2}},
\q\mbox{and}\q
X_0=0, ~I_{0-} =1,~ m_{0-}-\mbox{a.s.}
\eea
We next specify the function $a$,  which relies on  two functions $\k_i$ on $[0, T]\times \dbR$:
\bea
\label{eg-kappa}
\left.\ba{c}
\dis \k_0(t,x) := \dbE\big[\psi(x + W_t)\big],\q \k_1(t,x) := \dbE\big[\psi(x + W_t) \1_{\{W^*_t < 1\}}\big], \\
\dis a_t := {1\over 2}\big[ \k_0(t, 0) + t^2(1-t)^2\big]\1_{[0, 1]}(t) +  {1\over 2}\dbE\big[\k_1(t-1,  W_1)\big]  \1_{(1, 2]}(t),
\ea\right.
\eea
where $W^*_t := \sup_{0\le s\le t} W_s$. Recall Karatzas \& Shreve \cite[Chapter 2, Proposition 8.1]{KS} for the joint density of $(W_t, W^*_t)$, by direct calculations we have $0\le  \pa_t \k_0(t,x) - \pa_t \k_1(t,x)\to 0$ as $t\to 0$. Then
$
\pa_t \k_1(0,x)=\pa_t \k_0(0,x),
$
which implies that
$
a'_{1+} = h'_1 = a'_{1-},
$
that is, $a\in C^1([0, T])$.

\begin{prop}
\label{prop-eg-tau}
Under the above setting, an optimal $\dbP^*$ has the following structure:

\no {\rm (i)} At time $0$, there is a massive stop with  $\dbP^*(I_0=1)={1\over 2}$.

\no {\rm (ii)}There is no stop during the time interval $(0, 1]$: $I_t = I_0$, $0\le t\le 1$, $\dbP^*$-a.s.

\no {\rm (iii)} Particles stop continuously during the time interval $(1, 2)$:
\bea
\label{eg-tau}
\t = \inf\{t>1:  X_t - X_1 \ge 1\}\wedge 2 ,\q\dbP^*\mbox{-a.s. on}~\{I_0=1\}.
\eea

\no {\rm (iv)} All the remaining particles stop at time $2$.
\end{prop}
\proof  (i) Note that in this case 
\bea
\label{eg-a0}
v_0(m_{0-}) = \dbE^{m_{0-}}\big[\psi(X_0)I_{0-}\big]=\psi(0) =1 >   {1\over 2} = {1\over 2}\k_0(0, 0) = a_0.
\eea
Then $-(\dbL u_0+F)(0,m_{0-}) = v_0(m_{0-}) - a_0>0$, we have to stop some particles immediately. We may choose $m^*_0$ such that $m^*_0(I_0 = 1)={1\over 2}$, and then $t=0$ is a jump point of $m^*$, and
\beaa
v_0(m^*_{0}) =   \dbE^{m^*_{0}}\big[\psi(X_0)I_{0}\big]= {1\over 2},\q \mbox{and thus}\q -(\dbL u_0 + F)(0,m^*_{0}) = v_0(m^*_{0}) - a_0=0.
\eeaa  
Moreover, since $v_0(m^*_{0}) \in [a_0, v_0(m_{0-})]$, by \reff{eg-Cu0} we have $m^*_0\in C_{u_0}(0, m_{0-})$. This implies that $u_0(0, m^*_0) =u_0(0, m^*_{0-})$, and then it follows from \reff{optimality} that $\dbP^*$ is optimal at $t=0$.

(ii) For the $\dbP^*$ specified in the proposition, we have $I_t = I_0$ and hence $X_t = W^{\dbP^*}_t$ on $\{I_0=1\}$, $0\le t\le 1$, $\dbP^*$-a.s. By \reff{eg-kappa}  we see that
$v_0(m^*_t) = {1\over 2} \k_0(t, 0) \le a_t$, which implies that $-(\dbL u_0 + F)(0,m^*_t) = 0$, $0\le t\le 1$. Since no particle stops during this period, then by \reff{optimality} again $\dbP^*$ is optimal on $[0,1]$. 

(iii) We first note that, if we continue to keep all particles on $\{I_0=1\}$ alive after $t=1$, then we will have $v_0(m_t) = {1\over 2} \k_0(t, 0) > a_t$ (since $\k_0 > \k_1$) and thus  $-(\dbL u_0 + F)(0,m_t) >0$, which is not optimal. So after $t=1$, we start to stop particles, and our structure allows us to stop the particles continuously in the sense $m^*_t$ is continuous in $t$. Indeed, by \reff{eg-X} and \reff{eg-tau}, 
\beaa
\t = \inf\{t>1:  W^{\dbP^*}_t - W^{\dbP^*}_1 \ge 1\} \wedge 2 ,\q\dbP^*\mbox{-a.s. on}~\{I_0=1\}.
\eeaa
Then, for $t\in (1, 2)$,
\beaa
&&v_0(m^*_t)= \dbE^{\dbP^*}\Big[\psi(X_t) I_t\Big] = \dbE^{\dbP^*}\Big[\psi(X_t) I_0\1_{\{\t>t\}}\Big] \\
&&= \dbE^{\dbP^*}\Big[\psi(W^{\dbP^*}_1 + W^{\dbP^*}_t-W^{\dbP^*}_1) \1_{\{I_0=1\}} \1_{\{ \sup_{1\le s\le t}[W^{\dbP^*}_s - W^{\dbP^*}_1] < 1\}}\Big] ={1\over 2} \dbE^{\dbP^*}\Big[\k_1(t-1, W^{\dbP^*}_1) \Big] = a_t.
\eeaa
Therefore, $-(\dbL u_0 + F)(0,m^*_t) = 0$, $1< t<2$. 

Next, for any $1<t<2$, clearly $m^*_t$ is continuous, and thus $u_0(t, m^*_t) =u_0(t, m^*_{t-})$. Moreover, since $v_0(m^*_t)=a_t$, by \reff{eg-DIu0} and \reff{eg-phi} we have
\beaa
D_I u_0(t, m^*_t, X_t) =  [T-t]\pa_x   \f(t, v_0(m^*_t))\psi(X_t)  =[T-t]\pa_x   \f(t, a_t)\psi(X_t)  =0.
\eeaa
Then by \reff{optimality} again we see that $\dbP^*$ is optimal on $[1,2)$. 

(iv) This is required by our formulation of the problem.
\qed

\begin{rem}
\label{rem-eg-classical} {\rm
(i) For the $\dbP^*$ in Proposition \reff{prop-eg-tau}, $m^*_t$ has two jumps, one at $t=0$ and the other at $t=2$. In particular, the stopping at $t=0$ is randomized. Indeed, since $X_0\equiv 0$ under $m_{0-}$, there is no $A\in \cB(\dbR)$ such that $\dbE^{m_{0-}}[\psi(X_0)\1_A(X_0)] = a_0$.

\no(ii) If $X_0$ has continuous distribution under $m_{0-}$, say with density $\rho_0(x)$, then it is possible to have pure stopping strategy. Indeed, let $x_0$ be a median of $X_0$. Set $T$, $\psi$, $I_{0-}$, $\k_0$, $\k_1$ as in  \reff{eg-psi} and \reff{eg-kappa}, and modify the $a$ in \reff{eg-kappa} as follows:
\beaa
a_t := \Big[\int_{-\infty}^{x_0} \!\! \k_0(t, x) \rho_0(x) dx + t^2(1-t)^2\Big]\1_{[0, 1]}(t) + \int_{-\infty}^{x_0} \!\!\dbE[\k_1(t-1, x + W_1)] \rho_0(x) dx \1_{(1, 2]}(t).
\eeaa
By the same arguments as in Proposition \ref{prop-eg-tau}, the following pure stopping strategy is optimal:

$\bullet$ At time $0$, there is a massive stop for the particles $X_0 > x_0$: $I_0 = \1_{\{X_0\le x_0\}}$, $\dbP^*$-a.s.

$\bullet$ There is no stop during the time interval $(0, 1]$: $I_t = I_0$, $0\le t\le 1$, $\dbP^*$-a.s.

$\bullet$ Particles stop continuously during the time interval $(1, 2)$ following \reff{eg-tau}.

$\bullet$ All the remaining particles stop at time $2$. \qed}
\end{rem}

\section{Some extensions}
\label{sect-extension}

\subsection{Infinite horizon case}\label{infinite}

 This subsection is dedicated to the case $T = +\infty$. For any $(t,m) \in \ol \bQ_0$, let $\cP(t,m )$ denote the set of  $\dbP$ such that $\dbP_{Y_{t-}} = m$ and  \reff{asympt} holds on $[t, \infty)$.  We shall always assume
\begin{assum}
\label{assum-infinity}
{\rm (i)} Assumption \ref{assum-bsig} holds true on $[0, \infty)$;\\
{\rm (ii)}  $\int_0^\infty \sup_{m\in \cP_2(\bS)} |F(t, m)|dt <\infty$;\\
{\rm (iii)}  For any $(t,m)$ and $\dbP\in \cP(t,m)$, $X_\infty:= \lim_{t\to \infty} X_t$ exists, $\dbP$-a.s.
\end{assum}
We remark that one sufficient condition of (ii) above is that $|f(t,x, m)|\le Ce^{-\l t}$ for some constants $C, \l>0$, and a special case of (iii) is:
\bea
\label{GBM}
d=1,\q b = b_0 x, \q \si = \si_0 x,\q b_0 - {1\over 2}\si_0^2 <0.
\eea
That is,  the unstopped process $X^0$ in \reff{X0} is a Geometric Brownian motion and $X^0_\infty = 0$. 
 
We also define $I_\infty := 0$. This allows the case $\t = +\infty$, and guarantees that $\cP(t,m)$ is compact. The infinite horizon optimal stopping problem then simply writes:
\bea\label{optstop-inf}
V(t,m) := \sup_{\dbP \in \cP(t,m)}  \int_t^\infty F(s, \dbP_{Y_s}) ds + g(\dbP_{X_\infty}) \q \mbox{for all $(t,m) \in \ol \bQ_0$.}
\eea
The corresponding obstacle equation on Wasserstein space is
\bea\label{obstacle-inf}
\min_{m' \in C_u(t,m)}-[\dbL u + F](t,m') = 0, \ D_I u(t,m, \cdot) \ge 0, ~ (t,m) \in \ol \bQ_0,
\eea
with boundary condition $u(\infty,\cd) = g$.

Now by considering the problem on $[0, \infty]$, we see that all the definitions as well as all the results in the previous sections on the finite horizon remain true in the infinite horizon. 

\begin{rem}
\label{rem-homogeneous} {\rm
In the infinite horizon, one may naturally consider the time homogeneous case, that is, $b, \si, f$ do not depend on $t$. Then $V = V(m)$ is also time homogeneous, and thus \reff{obstacle-inf} becomes an elliptic problem: recalling \reff{DiffOp},
\bea\label{inf-elliptic}
\min_{m' \in C_u(m)} -\big[ \int_{\dbR^d}\!\!\!  \cL_x\d_m u_1(m',x) m'(dx,1) + F(m')\big] = 0, \ D_I u(m, \cdot) \ge 0, \ \mbox{for all $m \in \cP_2(\bS)$},
\eea
with boundary condition  $u = g$ on  $\pa \cP_2(\bS) := \{ m \in \cP_2(\bS) : m(\dbR^d, 1) = 0 \}$. We leave the details to interested readers. \qed}
 \end{rem}

\subsection{Mean field optimal stopping of a jump-diffusion}\label{JumpDiffSec}

This last subsection is dedicated to an informal discussion about the case where $(X,I)$ is a stopped jump-diffusion, i.e., $I_s = I_{t-} \1_{s < \t}$ and
\bea\label{jumpdiff}
X_s = X_t + \int_t^s b(r,X_r,m_r)I_r dr + \int_t^s \si(r,X_r,m_r)I_r dW_r + \int_t^s \g(r,X_{r-},m_{r-}) I_{r-} d\eta_r,
\eea
where $\eta$ is a pure jump process with intensity $\l_s := \l_s(s,X_s,m_s)$ and whose jump size is defined by a distribution $\n$, and $\g_s := \g(s,X_{s-},m_{s-})$ satisfies the usual conditions. We refer to Burzoni, Ignazio, Reppen \& Soner \cite{BIRS}, who characterized the mean field optimal control of a jump-diffusion by a dynamic programming equation (in the viscosity sense). The result of this section may be seen as a complement to the context of mean field optimal stopping. We consider the optimal stopping problem \eqref{weakoptstop}, where $\cP(t,m)$ is the set of probability measures such that the canonical process $(X,I)$ satisfies \eqref{jumpdiff}. Then, the value function still satisfies the DPP \eqref{weakDPP}. In order to formally derive the corresponding dynamic programming equation, we need to find the differential operator associated with the dynamics \eqref{jumpdiff}, which follows from Itô's formula \eqref{Ito0} in the present jump-diffusion case. Let $u \in C_2^{1,2}(\ol \bQ_0)$. Observing that the discontinuities in the flow $\bm=\{m_s\}$ are only due to $I$ (as $\eta$ has an intensity, hence no atoms), by shifting the jump at $s$ to $t$ as in \reff{JumpFunIto}, we have
\beaa
&&\dis u(s,m_{s-}) = u(t,m_{t^-}) + \int_t^s \dbL u(r,m_r)dr + \sum_{r \in J_{[t, s)}(\mathbf{m})} [u(r, m_r) -u(r, m_{r-})]  + \cJ_D,\\
&&\dis \mbox{where}\q  \cJ_D:=   \dbE^\dbP\Big[ \sum_{r \in J^c_{[t, s)}(\mathbf{m})}  \big(\d_m u(r, m_r, X_r, I_r)- \d_m u(t, m_r, X_{r-}, I_{r-})\big)\Big].
\eeaa

We next compute $\cJ_D$.  Denote  $\f_r(\cd) := \f(r,m_r,\cd)$ for any function $\f$ and $\D \eta_r:= \eta_r-\eta_{r-}$. Note again that $J_{[t, s)}(\mathbf{m})$ is countable and thus $\eta$ does not jump at   $J_{[t, s)}(\mathbf{m})$, a.s. Then
\begin{align*}
\cJ_D =& \ \dbE^\dbP\Big[\sum_{r \in J^c_{[t, s)}(\mathbf{m})}\big[ \d_m u_r(X_{r-} + \g_r(X_{r-})I_{r-} \D \eta_r,I_r) - \d_m u_r(X_{r-},I_{r^-}) \big]\Big] \\
=& \ \dbE^\dbP\Big[\sum_{r \in J^c_{[t, s)}(\mathbf{m})}  \big[ \d_m u_r(X_{r-} + \g_r(X_{r-})I_{r-} \D \eta_r,I_r) - \d_m u_r(X_{r-},I_{r}) \big] \Big] \\ 
 &+ \dbE^\dbP\Big[\sum_{r \in J^c_{[t, s)}(\mathbf{m})} \big[\d_m u_r(X_{r-},I_{r})   - \d_m u_r(X_{r-},I_{r^-})\big] \Big] \\
=& \ \dbE^\dbP\Big[\int_t^s \int_{\dbR^d} \big[\d_m u_r(X_{r} +y \g_r(X_{r})I_{r},I_r) - \d_m u_r(X_r, I_r)\big]\n(dy)\g_r(X_r) \l_r(X_r)I_r dr\ \Big] \\ &+ \dbE^\dbP\Big[ \int_{J^c_{[t, s)}(\mathbf{m})} D_I u_r(X_r)dI_r \Big],
\end{align*}
which implies that the differential operator corresponding to the dynamics \eqref{jumpdiff} is
\beaa
\dbL^{JD}u(t,m) := \dbL u(t,m) + \int_{(\dbR^d)^2}\!\!\!\! \big[ \d_m u_1(t,m,x + y\g(t,m,x))  - \d_m u_1(t,m,x)\big]\g \l(t,m,x) \n(dy)  m(dx,1).
\eeaa
Then, the dynamic programming equation corresponding to our problem is
\beaa
\underset{m' \in C_u(t,m)}{\min}\!\!\!\!-(\dbL^{JD} u + F)(t,m') = 0, 
~D_I u(t,m,.) \ge 0,~ u(T,\cd) = g,
~(t,m) \in [0,T] \!\times\! \cP_2(\mathbf{S}).
\eeaa
All the results of the previous sections can be adapted under appropriate assumptions. 

\appendix
\section{Proof of Proposition \ref{existence}}
\label{appendixA}

We assume for simplicity $t=0$ and fix $m\in  \cP_2(\bS)$. Let $C_m$ denote a generic constant which may depend on $T$ and $m$ but independent of $\dbP$. We proceed in three steps. 

{\it Step 1.} We first prove the following uniform integrability: denoting $X^*_T := \sup_{0\le s \le T} |X_s|$, 
\bea
\label{qUIlem}
\sup_{\dbP\in \cP(0,m)}\dbE^\dbP\big[ |X^*_T|^2\big]\le C_m,\qq \lim_{R \rightarrow \infty}\sup_{\dbP\in \cP(0,m)}\dbE^{\dbP}\big[ |X^*_T|^2 \1_{\{X^*_T \ge R\}}\big] = 0.
\eea
Indeed, for any $\dbP\in \cP(0, m)$, first by standard arguments we derive from \reff{asympt} that
$\dbE^\dbP\big[ |X^*_T|^2\big] \le C \dbE^m[1+|X_0|^2] \le C_m$.
In particular, this implies that the set $\{\dbP_{Y_s}: \dbP\in \cP(0, m), 0\le s\le T\}$ is bounded under $\cW_2$. Then, for any $p>2$, by  \reff{asympt}  again we have $\dbE^\dbP\big[ |X^*_T|^p\big|\cF_0\big] \le C_{m,p} [1+|X_0|^p]$, $\dbP$-a.s., where $C_{m,p}$ may depend on $p$ as well, but is still independent of $\dbP$. Now for any $R>0$,
\beaa
\dbE^\dbP\Big[|X^*_T|^2\1_{\{X^*_T\ge R\}}\Big] &&\le \dbE^\dbP\Big[|X^*_T|^2\1_{\{1+|X_0|\ge \sqrt{R}\}}\Big]+ \dbE^\dbP\Big[|X^*_T|^2\1_{\{{X^*_T\over 1+|X_0|}\ge \sqrt{R}\}}\Big]\\
&&\le \dbE^\dbP\Big[|X^*_T|^2\1_{\{1+|X_0|\ge \sqrt{R}\}}\Big]+ {1\over \sqrt{R}}\dbE^\dbP\Big[{|X^*_T|^3\over 1+|X_0|}\Big]\\
&&=\dbE^\dbP\Big[ \dbE^\dbP\big[|X^*_T|^2\big|\cF_0\big] \1_{\{1+|X_0|\ge \sqrt{R}\}}\Big]+ {1\over \sqrt{R}}\dbE^\dbP\Big[{\dbE^\dbP\big[|X^*_T|^3\big|\cF_0\big]\over 1+|X_0|}\Big]\\
&&\le C_{m,2}\dbE^m\Big[[1+|X_0|^2] \1_{\{1+|X_0|\ge \sqrt{R}\}}\Big]+ {C_{m,3}\over \sqrt{R}}\dbE^m\Big[1+|X_0|^2\Big].
\eeaa
Notice that the right side above does not depend on $\dbP$, then it clearly implies \reff{qUIlem}.

{\it Step 2.} We next show that $\cP(0,m)$ is closed under the weak convergence. Let $\{\dbP^n\}_{n \ge 1} \subset \cP(0,m)$ converge weakly to some $\dbP^\infty$. Since $\dbP^n_{(X_0,I_{0-})} = m$ for all $n$, we have  $\dbP^\infty_{Y_{0-}} = m$. Then it suffices to show that the processes $M,N$ in \eqref{martingalepb} are $\dbP^\infty-$martingales on $[0,T]$. We shall report only the detailed argument for $M$, as it is immediately adapted to $N$. 

Notice that the support of $\dbP^\infty$ is separable under the Skorokhod distance $d_{SK}$, as a subspace of the separable metric  space $\O$. Then it follows from the Skorokhod's representation theorem, see Billingsley \cite[Theorem 6.7]{Bil}, that there exists a probability space $(\O^0, \cF^0, \dbP^0)$ and processes $\{Y^n := (X^n,I^n)\}_{n \ge 1}$ and $Y^\infty := (X^\infty, I^\infty)$ defined on this space such that,  
 \bea
 \label{XInconv}
 \dbP^n_{Y} = \dbP^0_{Y^n} \ \mbox{for all $n \le \infty$, and} \ d_{SK}(Y^n, Y^\infty) \underset{n \rightarrow \infty}{\longrightarrow} 0, \ \dbP^0-\mbox{a.s.}
 \eea
 For all $n \ge 1$, the $\dbP^n-$martingale property of $M$ translates to:
 \bea\label{projsko}
  \dbE^{\dbP^0}[(M_s^n-M_{t}^n)\psi(Y_{. \wedge {t}}^n)] = 0 \q \mbox{ for all $\psi \in C_b(\O)$ and $0\le t\le s\le T$,}
 \eea
 with $M^n_s = X^n_{t} - \int_{t}^s b(r,X^n_r,\dbP^0_{Y^n_r})I_r^ndr$ and $C_b(\O)$ the set of $\dbR^d$-valued bounded continuous functions on $\O$. Moreover,  for $r \in [t,T]$, by the Lipschitz continuity of $b$ we have
 \beaa
 \lvert b(r,X_r^n, \dbP^0_{Y_r^n}) - b(r,X_r^\infty, \dbP^0_{Y_r^\infty})\rvert &\le& C\big[\lvert X_r^n - X_r^\infty \rvert + \cW_2(\dbP^0_{Y_r^n},\dbP^0_{Y_r^\infty})\big].
 \eeaa
Send $n\to \infty$, by \reff{XInconv} we have $\lvert X_r^n - X_r^\infty \rvert \to 0$, $\dbP^0$-a.s. and  $\dbP^0_{Y_r^n} \to \dbP^0_{Y_r^\infty}$ weakly. Then, by the 2-uniform integrability \reff{qUIlem} of $\{\dbP^0_{Y_r^n}\}_{n \ge 1}$, we have
  $\cW_2(\dbP^0_{Y_r^n},\dbP^0_{Y_r^\infty}) \to 0,$ see Carmona \& Delarue \cite[Vol. I, Theorem 5.5]{CarDel}.
  Thus
  $$  b(r,X_r^n, \dbP^0_{Y_r^n}) \underset{n \rightarrow \infty}{\longrightarrow} b(r,X_r^\infty, \dbP^0_{Y_r^\infty}), \ \mbox{$\dbP^0$-a.s.} $$
  Moreover, as $b$ is Lipschitz and $\{X^n\}_{n \ge 1}$ are uniformly integrable (as the 2-uniform integrability of \reff{qUIlem} implies the 1-uniform integrability), then $\{M^n\}_{n\ge 1}$ are uniformly integrable. The convergence for the Skorokhod distance also implies the convergence of $I_{. \wedge t}^n$ to $I_{. \wedge t}^\infty$.
This allows to take the limit in \eqref{projsko} as $\psi \in C_b(\O)$, hence
 $ \dbE^{\dbP^0}[(M_s^\infty-M_t^\infty)\psi(Y_{. \wedge t}^\infty)] = 0.$
By the arbitrariness of $\psi \in C_b(\O)$, this proves $M^\infty$ is a $\dbP^0-$martingale, or equivalently that $M$ is a $\dbP^\infty-$martingale. 

{\it Step 3.} We now show that $\cP(0, m)$ is compact under $\cW_2$. Let $\{\dbP^n\}_{n\ge 1}\subset \cP(0, m)$. First, by the first estimate in \reff{qUIlem} and noticing that $I$ is bounded by $1$, one can easily obtain a uniform bound for the conditional variation of $Y$ under all $\dbP^n$, then by Meyer \& Zheng \cite[Theorem 4]{ZM} we see that $\{\dbP^n\}_{n\ge 1}$ is relatively weakly compact, namely there exists a weakly convergent subsequence. By Step 2,  without loss of generality we assume the whole sequence $\dbP^n\to \dbP^\infty\in \cP(0, m)$ weakly. Moreover, by the second estimate in \reff{qUIlem} $\{\dbP^n\}_{n\ge 1}$ is 2-uniformly integrable, then it follows from Carmona \& Delarue \cite[Vol. I, Theorem 5.5]{CarDel} again that $\underset{n \rightarrow \infty}{\lim}\cW_2(\dbP^n,\dbP^\infty) = 0$. This proves the compactness of $\cP(0, m)$.

Finally, since $g$ is upper-semicontinuous, the above compactness implies the existence of optimal $\dbP^*$ for the mean field optimal stopping problem \reff{weakoptstop}.
 \qed

\section{Proof of Theorem \ref{Itothm}}\label{appendixB}

Let $\Xi_\bm$ denote the convex hull of $\{m_s, m_{s-}: 0\le s\le T\}$: 
\beaa
\Xi_\bm := \Big\{\l m_{s'} + (1-\l) m_{t'}: 0\le \l\le 1, 0\le s\le t\le T, s'=s, s-, t' = t, t-\Big\} \subset \cP_2(\dbR^{d'}).
\eeaa  
We first show that $\Xi_\bm$ is compact. Indeed, for any $(\l_n, s_n', t_n')$, there exists a convergent subsequence and we may assume without loss of generality that $(\l_n, s_n, t_n) \to (\l, s, t)$. By considering different cases, one can easily show that, possibly along a subsequence, for some $s', t'$ we have $\l_n m_{s'_n} + (1-\l) m_{t'_n}\to \l m_{s'} + (1-\l) m_{t'}\in \Xi_\bm$, thus $\Xi_\bm$ is compact.

Denote $\D Y_s:= Y_s-Y_{s-}$ and $Y^D_t := \sum_{0<s\le t} \D Y_s$. By \reff{summability} it is clear that
\bea
\label{EY*}
\dbE\big[|Y^*_T|^2 + \|Y^D\|_T^2\big]<\infty,\q\mbox{where}\q Y^*_T := \sup_{0\le s\le T} |Y_s|,\q  \|Y^D\|_t := \sum_{0<s\le t} |\D Y_s|.
\eea
 For $n\ge 1$, set $\D t := {T\over n}$, $t_i:= i \D t$, $i=0,\cds, n$. Then, for each $i$, 
\bea
\label{Du}
\left.\ba{c}
\dis u(t_{i+1}, m_{t_{i+1}})-u(t_i, m_{t_i}) =\int_{t_i}^{t_{i+1}} \pa_t u(s, m_{t_{i+1}})ds + \int_0^1 \dbE[\xi^\l_{t_{i+1}}] d\l,\\
\dis \mbox{where}\q \xi^\l_{t_{i+1}}:= \d_m u(t_i, m_{t_i}^\l,Y_{t_{i+1}})-\d_m u(t_i, m_{t_i}^\l,Y_{t_i}),\q m^\l_{t_i} := \l m_{t_i} + [1-\l]m_{t_{i+1}}.
\ea\right.
\eea
By the standard It\^{o}'s formula:
\beaa
&\dis \xi^\l_{t_{i+1}} = \int_{t_n}^{t_{n+1}}  \big[\G^{2,\l}_s\cd dY_s^c  +  \frac{1}{2} \G^{3,\l}_s:d\langle Y^c \rangle_s\big] + \int_{(t_n, t_{n+1}]} \G^{4,\l}_s dY^D_s,\\
&\dis \mbox{where} \q \G^1_s :=   \pa_t u(s, m_{t_{i+1}}),\q  \G^{2,\l}_s :=  \pa_y \d_m u(t_i, m_{t_i}^\l,Y_s),\q \G^{3,\l}_s :=  \pa_{yy}^2 \d_m u(t_i, m_{t_i}^\l,Y_s),\\
&\dis \G^{4,\l}_s:= \int_0^1 \pa_y \d_m u\big(t_i, m_{t_i}^\l, \th Y_s + [1-\th] Y_{s-}\big) d\th.
\eeaa
Note that $m_{t_{i+1}}, m^\l_{t_i}\in \Xi_\bm$, by the growth conditions in Definition \ref{linderiv} we have 
\bea
\label{Gammabound}
 |\G^1_s|\le C,\q |\G^{2,\l}_s| \le C[1+|Y_s|],\q  |\G^{3,\l}_s|\le C,\q  |\G^{4,\l}_s|  \le C[1+|Y_s|+|Y_{s-}|].
\eea
Then
\beaa
\dbE\Big[ \Big(\int_{t_n}^{t_{n+1}}\G^{2,\l}_s(\G^{2,\l}_s)^\top  : d\la M^c\ra_s\Big)^{1\over 2}\Big] \le C\dbE\Big[ \Big([1+|Y^*_T|^2] \la M^c\ra_T\Big)^{1\over 2}\Big] \le C\dbE\Big[ 1+|Y^*_T|^2 +\la M^c\ra_T\Big]<\infty.
\eeaa
This implies $\int_0^1 \dbE\big[\int_{t_n}^{t_{n+1}} \G^{2,\l}_s \cd dM^c_s\big]d\l=0$, and thus
\beaa
u(T, m_T) = u(0, m_0) +  \int_0^T \G^1_s ds +  \int_0^1\dbE\Big[\int_0^T \big[\G^{2,\l}_s \cd d A^c_s + \G^{3,\l}_s : d\la Y^c\ra_s\big] + \int_{(0, T]} \G^{4,\l}_s dY^D_s\Big] d\l.
\eeaa
Fix $\l$, $s$, and send $n\to \infty$. By the regularity of $u$ we have:  denoting $m^\l_s:= \l m_{s-} + [1-\l] m_s$,
\beaa
&\G^1_s \to \pa_t u(s, m_s),\q \G^{2,\l}_s\to   \pa_y \d_m u(s, m^\l_s,Y_s),\q \G^{3,\l}_s \to  \pa_{yy}^2 \d_m u(s, m^\l_s,Y_s),\\
& \G^{4,\l}_s\to \int_0^1 \pa_y \d_m u\big(s, m_s^\l, \th Y_s + [1-\th] Y_{s-}\big) d\th,\q a.s.
 \eeaa
By \reff{summability}, \reff{EY*}, and \reff{Gammabound},  we may apply the dominated convergence theorem to obtain
\beaa
&&u(T, m_T) = u(0, m_0) +  \int_0^T \pa_t u(s, m_s)  ds +  \int_0^1\dbE\Big[\int_0^T \big[ \pa_y \d_m u(s, m^\l_s,Y_s) \cd d A^c_s \\
&&+  \pa_{yy}^2 \d_m u(s, m^\l_s,Y_s) : d\la Y^c\ra_s\big] + \int_{(0, T]} \int_0^1 \pa_y \d_m u\big(s, m_s^\l, \th Y_s + [1-\th] Y_{s-}\big) d\th dY^D_s\Big] d\l.
\eeaa
Since $A^c$ and $\la Y^c\ra$ are continuous, and $m_s$ has at most countably many jumps, then
\bea
\label{Ito1}
&&u(T, m_T) = u(0, m_0) +   \int_0^T \pa_t u(s, m_s)  ds +  \cJ_D  \nonumber\\
&&\q +\dbE\Big[\int_0^T\!\!\! \big[ \pa_y \d_m u(s, m_s,Y_s) \cd d A^c_s+  \pa_{yy}^2 \d_m u(s, m_s,Y_s) : d\la Y^c\ra_s\big]\Big], \\
&&\mbox{where}\q \cJ_D:=  \dbE\Big[ \int_0^1 \int_{(0, T]} \int_0^1 \pa_y \d_m u\big(s, m_s^\l, \th Y_s + [1-\th] Y_{s-}\big) d\th dY^D_s d\l\Big] .\nonumber
\eea

It remains to compute $\cJ_D$. First, by Fubini's theorem,
\bea
\label{JD0}
&\dis\cJ_D =  \dbE\Big[ \sum_{s\in (0, T]}\int_0^1  \int_0^1 \pa_y \d_m u\big(s, m_s^\l, \th Y_s + [1-\th] Y_{s-}\big) d\th  d\l \D Y_s\Big] =\dbE\Big[ \sum_{s\in (0, T]} \D \d_m u_s\Big]\nonumber\\
&\dis\mbox{where}\q \D \d_m u_s:= \int_0^1 \big[ \d_m u(s, m_s^\l,  Y_s ) - \d_m u(s, m_s^\l,  Y_{s-} )\big]  d\l.
\eea
Note that $(0, T] = J_{(0, T]}(\bm) \cup J^c_{(0, T]}(\bm)$. Since $J_{(0, T]}(\bm)$ is countable, then
\bea
\label{JD1}
\dbE\Big[ \!\!\! \sum_{s\in J_{(0, T]}(\bm) } \!\!\! \D \d_m u_s\Big]=  \!\!\!\sum_{s\in J_{(0, T]}(\bm) } \!\!\!\dbE\big[ \D \d_m u_s\big] =  \!\!\!\sum_{s\in J_{(0, T]}(\bm) } \!\!\! \big[\d_m u(s, m_s) - \d_mu(s, m_{s-})\big],
\eea
where the second equality is due to \reff{dmu}. Next, for $s\in J^c_{(0, T]}(\bm)$, we have $m^\l_s = m_s$, $0\le \l\le 1$. Then
\bea
\label{JD2}
\dbE\Big[ \!\!\! \sum_{s\in J^c_{(0, T]}(\bm) } \!\!\! \D \d_m u_s\Big]= \dbE\Big[ \!\!\! \sum_{s\in J^c_{(0, T]}(\bm) } \!\!\! \big[ \d_m u(s, m_s,  Y_s ) - \d_m u(s, m_s,  Y_{s-} )\big]\Big].
\eea
We emphasize that, since $J^c_{(0, T]}(\bm)$ is uncountable, unlike in \reff{JD1} we cannot switch the order of $\dbE$ and $\sum_{s\in J^c_{(0, T]}(\bm) }$ at above. Now plug \reff{JD1}, \reff{JD2} into \reff{JD0}, and then plug \reff{JD0} into \reff{Ito1}, we complete the proof. \qed

\bibliographystyle{plain}
\bibliography{references}
\end{document}